# Optimized dynamic scheduling of an exclusive bus lane or high occupancy vehicle lane in a bimodal traffic corridor

Jiaqing Lu[a], Dian Sheng[b], Qianwen (Vivian) Guo[a*] and Paul Schonfeld[c]

[a]*Department of Civil and Environmental Engineering, Florida State University, Tallahassee, United States;*

[b]*School of Management, Huazhong University of Science and Technology, Wuhan, China;*

[c]*Department of Civil and Environmental Engineering, University of Maryland, College Park, United States*

*corresponding author: Qianwen (Vivian) Guo, Department of Civil and Environmental Engineering, Florida State University, Tallahassee, United States, qguo@fsu.edu

# Optimized dynamic scheduling of an exclusive bus lane or high occupancy vehicle lane in a bimodal traffic corridor

## Abstract

Efficient management of traffic corridors is critical for sustaining urban mobility. Exclusive bus lane (EBL) and high occupancy vehicle lane (HOVL) are two prominent strategies for enhancing public transit services and alleviating congestion. EBLs prioritize bus transit by providing dedicated lanes for faster travel times, while HOVLs encourage carpooling by reserving lanes for high-occupancy vehicles. However, static implementations of these policies may underutilize road resources and disrupt general-purpose lanes. Dynamic implementation, based on real-time demand, can potentially maximize road efficiency and minimize negative impacts. This study first compares the mixed traffic policy (MTP), exclusive bus lane policy (EBLP), and high occupancy vehicle lane policy (HOVLP) by formulating the total system costs in the context of a bimodal traffic corridor involving private cars and public buses. Under each lane policy, the bus frequency is optimized together with the modal split equilibrium derived separately. Based on dynamic demand simulated using an Ornstein-Uhlenbeck (O-U) process, switching thresholds are then derived to identify optimal periods for implementing each policy. Results reveal significant reductions in total system costs with the proposed dynamic policy schedules. Compared to static implementations, the dynamic policy schedules achieve cost reductions of 11.9%, 6.70%, and 43.64% relative to MTP-only, EBLP-only, and HOVLP-only scenarios, respectively. Additionally, in two real case studies of existing EBL and HOVL operations in Seattle, the proposed dynamic policy reduces total costs by 32.5% and 28.6%, respectively. The findings provide valuable insights for policymakers and transit planners, offering a robust framework for dynamically scheduling and integrating EBL and HOVL policies to optimize urban corridor efficiency and reduce overall system costs.

## 1. Introduction

Urban transportation corridors face mounting challenges as rising travel demand collides with limited roadway capacity (Li et al., 2012; Wang and Du, 2016). The coexistence of private automobiles and public transit creates direct competition for scarce road space, leading to congestion that undermines operational efficiency and service quality for all users (Yao et al., 2012 and 2015). To mitigate these issues, city administrators implement various traffic management policies to optimize road infrastructure usage and strengthen public transit (Cats, 2011; Li et al., 2020). Among others, two prominent strategies are Exclusive Bus Lanes (EBL) and High Occupancy Lanes (HOVL). EBLs, introduced in U.S. cities during the 1960s and 1970s, prioritize bus transit by dedicating specific lanes for exclusive use, thereby aiming to improve operational speeds and schedule reliability (Basso et al., 2011). HOVLs, first established in 1969, encourage carpooling and the use of other high-occupancy modes by reserving lanes for vehicles that meet minimum occupancy thresholds, thereby discouraging single-occupancy travel (Cohen et al., 2022). While EBLs directly prioritize public transit, HOVLs aim to shift travel behavior toward higher-occupancy modes (Jang et al., 2014). Both strategies, however, require the reallocation of limited roadway space, a decision that inevitably entails trade-offs and potential conflicts among different road user groups (Huang et al., 2025).

The intended benefits of these policies are significant. EBLs can reduce bus travel times and passenger waiting times, improve service reliability and passenger satisfaction, reduce accident rates and emissions, and attract more individuals to public transit (Kim et al., 2019; Zhang, 2023; Li et al., 2025). For instance, Seattle's 3rd Avenue EBL has notably reduced delays, with studies documenting measurable gains in bus service efficiency (Bancroft, 2024). HOVLs aim to reduce SOV usage, thereby improving travel time predictability and reliability for carpoolers, offering environmental benefits through emission reductions, and fostering a culture of ridesharing by providing a reduced congestion ride incentive (Dahlgren et al., 1998; Brownstone et al., 2003; Boriboonsomsin and Barth, 2008; Lou, 2013; Fontes et al., 2014; Javid et al., 2017; Wiseman, 2019; Guo and Ma, 2020; Zhong et al., 2020; Boysen et al., 2021;

Cohen et al., 2022). A critical, often understated, benefit for both is the enhancement of travel time reliability, which is a major factor in commuter satisfaction and mode choice (Yao et al., 2012; Wang et al., 2019). Ultimately, the long-term effectiveness of these strategies depends on their ability to induce a significant modal shift, which is key to achieving sustainable, long-term congestion relief.

Despite these advantages, static implementations of EBL and HOVL policies, e.g., operating permanently or on fixed schedules, face important limitations. A key concern is the potential for underutilization of dedicated road resources, particularly during off-peak periods or when demand for the prioritized mode is low (Khoo et al., 2014; Wu et al., 2017). This empty lane syndrome, where a dedicated lane remains largely unused while adjacent general-purpose lanes are heavily congested, can erode public and political support for such measures. Furthermore, static policies may worsen congestion in general-purpose lanes, leading to concerns about the equitable distribution of road resources and potentially negating overall system benefits (Chu et al., 2012; Yao et al., 2012). Indeed, some EBL projects have been discontinued due to worsening general traffic conditions, and certain studies suggest that HOVLs, when active during peak congestion, might decrease overall traffic capacity compared to unrestricted lane usage (Menendez and Daganzo, 2007). This highlights a critical balance, whereby a poorly calibrated static policy can lead to loss of system efficiency if the negative impacts on general traffic outweigh the gains on the dedicated lane.

Another challenge is that urban travel demand is inherently dynamic, fluctuating by time of day, day of week, and in response to incidents or special events. Static policies lack the flexibility to respond to these variations, often leading to periods of inefficiency (Othman et al, 2023). For example, current guidelines for EBL implementation often rely on generalized peak traffic patterns (Othman et al, 2025). Thus, the inherent limitations of static approaches underscore the potential of dynamic control strategies for EBL and HOVL management. By adjusting lane allocations and operational rules in real time, dynamic control offers the promise of maximizing network efficiency (Zhao and Zhou, 2018). This involves leveraging Intelligent Transportation Systems (ITS), including sensors, data analytics, and Variable Message Signs (VMS), to adaptively adjust lane allocations or operational rules (Zhao et al, 2019).

Dynamic control also enhances resilience by reallocating roadway capacity during unexpected events such as accidents (Lai et al., 2020). One notable application is the intermittent EBL, where lanes operate exclusively for buses only under specific conditions. Typically, when a bus enters or approaches the lane, other vehicles are required to merge out, freeing the lane. Zhao and Zhou (2018) show that intermittent exclusive bus lanes (EBLs) can reduce average intersection delays by over 10%. In particular, integrating connected vehicles and V2X communication with intermittent-priority bus lanes further improves bus reliability and efficiency while minimizing adverse effects on general traffic across different demand scenarios (Kampitakis et al., 2024). Szarata et al. (2021) demonstrate that dynamic bus lanes yield better bus service benefits with minimal private car delays compared to permanent EBLs. Successful implementations of intermittent EBLs in cities such as Lisbon, Melbourne, and Eugene have improved transit speeds by 8-20% (Viegas et al., 2006; Currie and Lai, 2008). However, intermittent lane models face practical implementation challenges, which often rely on preset schedules based on typical traffic patterns, such as peak and off-peak hours, rather than on quantitatively optimized operational periods (Chiabaut & Barcet, 2019).

Despite promising evidence, significant research gaps remain. Many existing studies focus on network-wide or intersection-based applications, often with static settings, or prescribe policies based on general peak patterns rather than quantitatively derived optimal operational periods for specific corridors (Othman et al., 2023). A critical deficiency is the lack of an integrated decision framework that cohesively considers the complex interdependencies among Mixed Traffic Policy (MTP), EBL Policy (EBLP), and HOVL Policy (HOVLP), as depicted in Figure 1. Key operational parameters, such as bus frequency, demand split under user equilibrium, and stochastic variations in travel demand, are rarely considered.

Foremost among these parameters is bus frequency, whose optimization within each specific lane policy is fundamental to achieving desired outcomes. The choice of lane policy directly influences bus performance, and optimal frequency must therefore be tailored to the active policy. For example, an EBLP might enhance bus speeds and attract more riders, requiring higher frequency to maintain service quality and fully realize its benefit. Conversely, under an MTP, a lower frequency may balance service attractiveness with operational efficiency in mixed traffic. Existing analytical frameworks frequently overlook these crucial dynamic interactions, including how bus frequency should be adjusted in response to fluctuating demand under different lane policies, alongside optimized demand split for bus and auto according to the inherent stochasticity of travel demand.

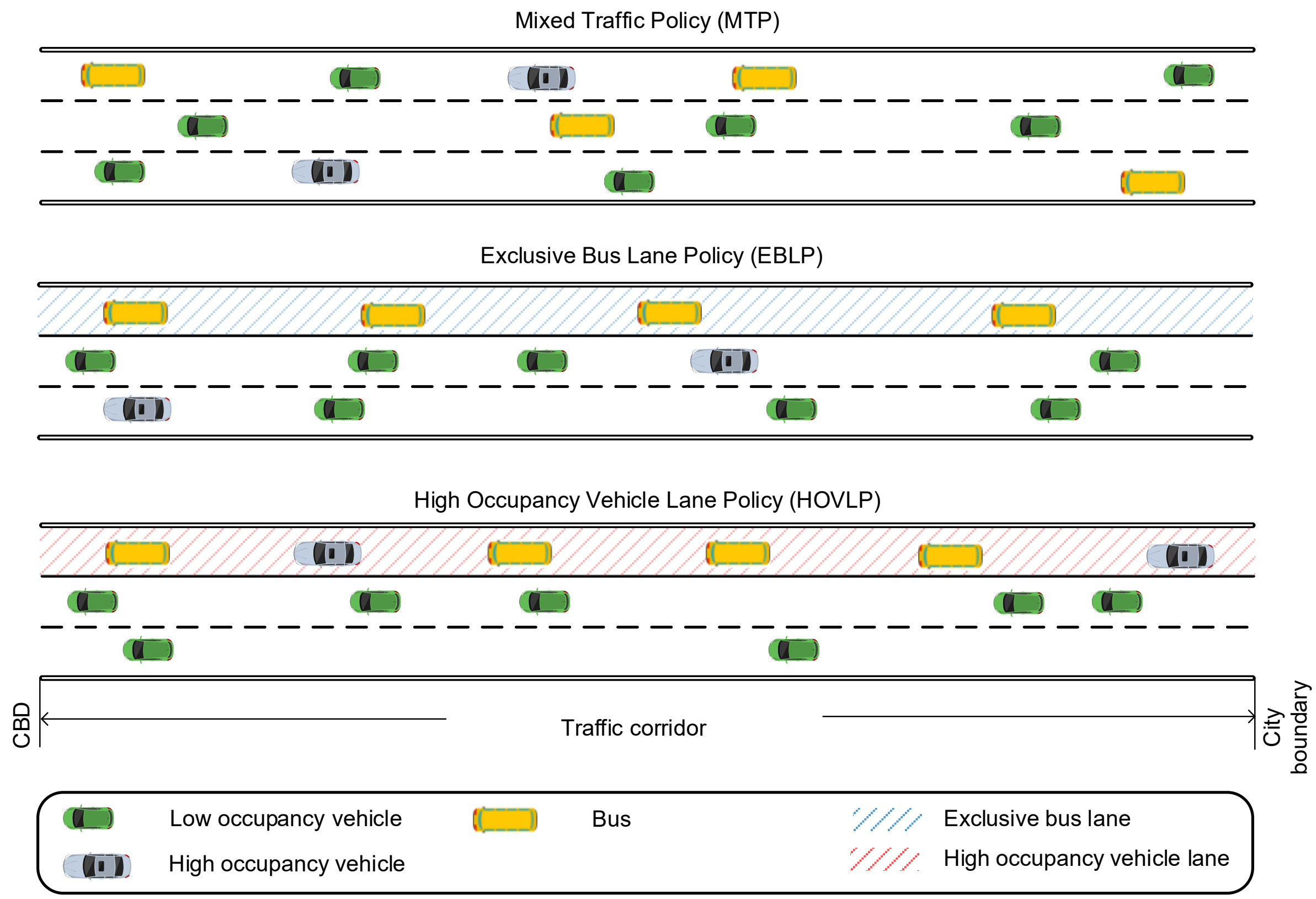

**Figure 1**. Illustration of different lane management policies in a traffic corridor.

Moreover, uncertainty in travel demand further complicates policy design. Corridor demand fluctuates not only with predictable daily patterns but also due to random factors such as incidents, delays, or traveler behavior. Such stochastic variations, if ignored, limit the effectiveness of deterministic models. To address this, we model travel demand using the Ornstein-Uhlenbeck (O-U) process, which effectively captures mean-reverting stochastic behavior over time (Berling & Martínez-de-Albéniz, 2011; Guo et al., 2017). This stochastic modeling approach enables us to develop robust switching thresholds for policy transitions. Nevertheless, the critical question of when to switch between different lane policies and how to explicitly incorporate the optimized bus frequency for each potential policy to guide these transitions remains largely unresolved. Most prior studies have relied on inflexible and predetermined thresholds, which fail to integrate the crucial element of policy-specific bus frequency optimization into the dynamic switching decision (Othman et al., 2023; Yao et al., 2018; Yang, 1998; Zhou et al., 2020).

To address these gaps, this study develop an analytical framework for the dynamic switching and scheduling of MTP, EBLP, and HOVLP along a bi-modal transportation corridor, in which commuters can choose between auto and bus. The framework integrates cost functions, endogenous bus frequency optimization, and stochastic demand modeling to support adaptive policy decisions. The major contributions are fourfold:

(i) We formulate enhanced cost functions for evaluating MTP, EBLP, and HOVLP, incorporating passenger discomfort and policy implementation costs in addition to conventional travel and waiting costs;

(ii) We propose a system cost minimization model that endogenously optimizes bus frequency, accounting for trade-offs between user costs and operator costs;

(iii) We develop a dynamic scheduling model that incorporates stochastic demand fluctuations, modeled via the O-U process, to capture realistic temporal variability;

(iv) We design three distinct policy schedule frameworks and demonstrate that dynamic management significantly outperforms static approaches in reducing costs and improving efficiency.

The rest of this paper is organized as follows. Section 2 describes the travel demand and formulates the system cost minimization problem. Section 3 introduces the policy switching and scheduling model for optimizing the timing and duration of different lane policies. Section 4 presents results from simulations and case studies, highlighting the performance gains of dynamic scheduling. Finally, Section 5 concludes and discusses directions for future research.

## 2. Problem description

In this section, the objective is to specify the total system costs under three lane management policies and optimize bus frequency statically. In the next section, a scheduling model that dynamically switches between alternative lane policies will be proposed to exploit the advantage of each policy under varying travel demand.

Consider a traffic corridor connecting the city's Central Business District (CBD) and its suburban area, represented as a linear segment, a common approach in transportation studies (Liu et al., 2009; Sun et al., 2017; Guo et al., 2021; Guo et al., 2023; Yu et al., 2023). As illustrated in Figure 2, the bimodal corridor includes two alternative travel modes: private automobiles and buses. The commute demand within this corridor is characterized by a many-to-one pattern, with multiple commuters heading towards the CBD (Zhang et al., 2005). It is assumed that commuters select the most cost-effective transportation mode (i.e., bus or auto) for their journey to the CBD. The government's role is to design and implement lane management policies to optimize corridor performance. This involves deciding which lane policy (MTP, EBLP, HOVLP) to activate and, crucially, determining the optimal bus frequency for the public transit agency operating the bus services. The government must then implement the chosen policy, potentially using ITS infrastructure, and monitor its effectiveness, making adjustments as necessary.

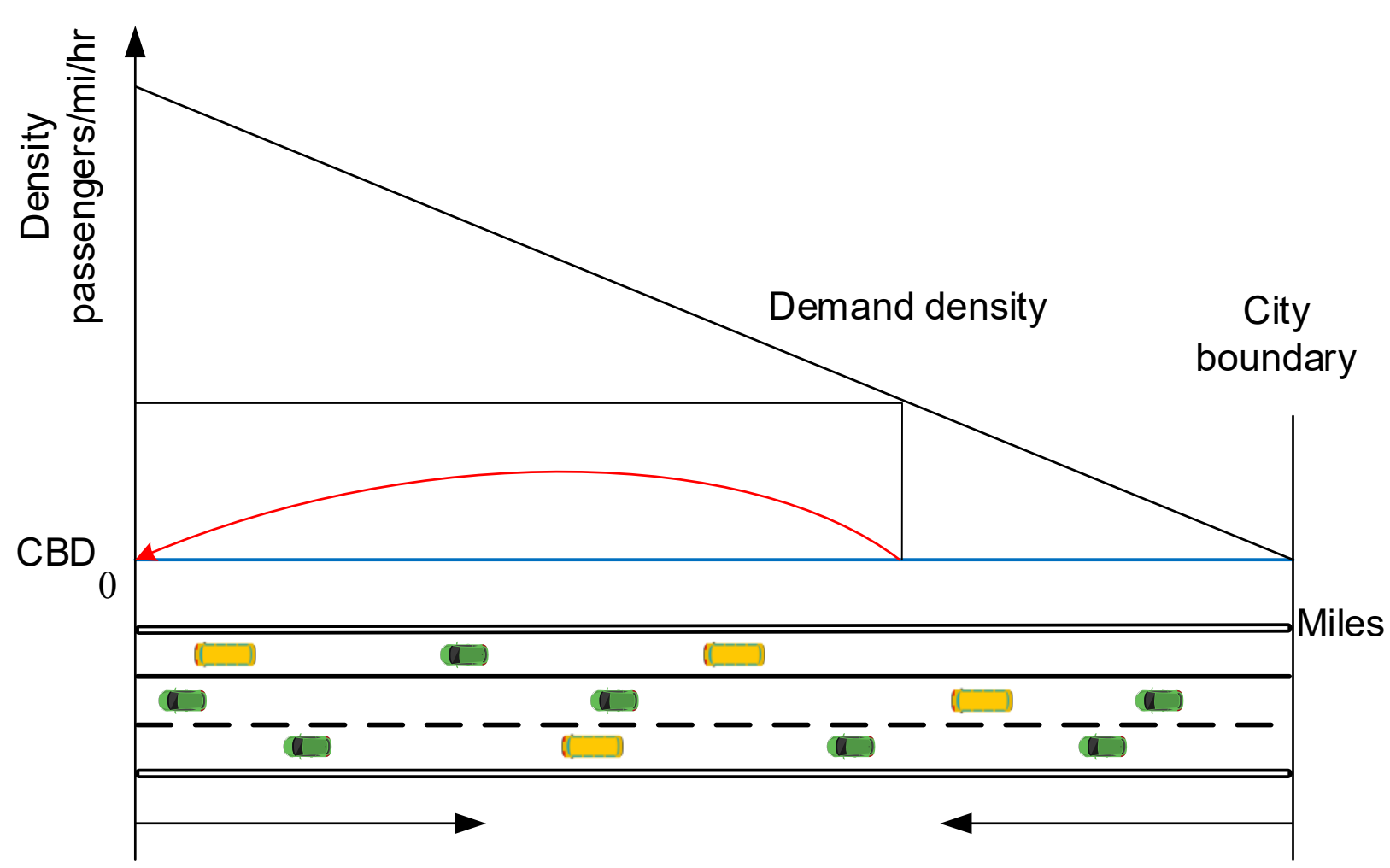

**Figure 2**. Illustration of a bimodal traffic corridor.

In the benchmark scenario, the auto and bus share the same road lanes under the MTP. To enhance traffic efficiency, the government or relevant public agency considers two alternative lane management policies: the EBLP, which dedicates one lane exclusively to buses, and the HOVLP, which reserves a lane for vehicles with multiple occupants, including buses. The bus services within the corridor are assumed to be operated by a public transit agency. To evaluate the system cost of each alternative, we first formulate the commuters' travel cost for auto and bus under different policies. A deterministic continuum mode choice equilibrium model is then derived to characterize commuters' mode choices. Based on that, one can calculate the system cost of implementing each policy, which includes all commuters' travel costs, bus operating costs, and policy implementation costs. This static evaluation provides the basis for optimizing bus frequency and comparing the overall performance of the three strategies under steady-state conditions.

All notations used in the formulation of our model are clearly defined in Table 1.

**Table 1**. Notation list.

| Parameters | Definitions |
|---|---|
| $A$ | Distance between CBD and city boundary |
| $C$ | Lane capacity |
| $C_O^f, C_O^v$ | Fixed/ Variable cost of bus operation |
| $C_e^f, C_e^v, C_h^f, C_h^v$ | Fixed/ Variable cost of EBL/HOVL signal implementation |
| $C_p^a, C_p^b$ | Total travel cost of auto/bus users under different policies |
| $C_p^O$ | Bus operating cost under different policies |
| $C_p^s$ | EBL/HOVL signal implementation cost |
| $F_p$ | Bus service frequency under different policies |
| $K$ | Factor converting a bus into a number of autos with equivalent effect on traffic |
| $f_b$ | Flat fare paid per bus trip |
| $n_{lane}$ | The number of lanes in the traffic corridor |
| $n_\Gamma$ | The number of schedule periods in the set $\Gamma$ |
| $O_a$ | Average auto occupancy |
| $O_{la}, O_{ha}$ | Occupancy rate of low-occupancy/ high-occupancy auto |
| $O_b$ | Bus capacity |
| $q_0$ | Demand density in the CBD |
| $t^{entry}, t^{exit}$ | Starting/ Exiting point from one policy to another |
| $t_0^a, t_0^b$ | Free-flow auto/bus travel time per unit of distance |
| $W$ | Cumulative cost in the schedules period |
| $Z$ | Total cost saving |
| $\alpha^a, \alpha^b, \beta^a, \beta^b$ | Corrected parameter of BPR function |
| $\varepsilon_t^a, \varepsilon_t^b$ | Value of in-vehicle time of auto users/bus travelers |
| $\delta_f^a, \delta_v^a$ | Fixed cost/ Variable cost of auto usage |
| $\iota_1, \iota_2$ | Positive parameters in discomfort function |
| $\mu$ | The fraction of low-occupancy autos over total autos |
| $\gamma_1, \gamma_2, \gamma_3$ | Positive calibrated parameters in waiting time function |
| $\varrho_l, \varrho_h$ | Fraction of low-occupancy/ high-occupancy auto travelers over total number of auto travelers |
| $\upsilon$ | Mean reversion parameter |
| $\varpi$ | Long-term ratio of the O-U process |
| $\sigma$ | Volatility of the O-U process |
| $\Gamma$ | Feasible lane policy schedule set |

### *2.1 Travel demand*

As shown in Figure 2, the density of travel demand at location *x* is:

$$q(x) = q_0 - \frac{q_0}{A}x, \ \forall \, x \in (0, A], \tag{1}$$

$q_0$ is the demand density in the CBD, and *A* is the length of the corridor, with demand density becoming zero at the city boundary (Li et al., 2015; Sun et al., 2017; Guo et al., 2018). During the morning peak hours, all commuters choose the less costly mode when traveling along the corridor to CBD. When everyone selects the least costly mode, the congested bimodal corridor reaches an equilibrium state (Liu et al., 2025). In such a state, the generalized travel cost of any mode at a given location must be the lowest among the two modes if commuters at that location choose to use it. Therefore, the equilibrium mode choice is used to split the demand in the bimodal corridor represented by Eq. (2). Appendix A provides the formulation and proof of this deterministic continuum equilibrium mode choice.

$$q(x) = {q_a}^*(x) + {q_b}^*(x) = \ R^* q(x) + (1 - R^*) q(x), \forall \, x \ \in (0, A]. \tag{2}$$

where $q_a^*(x)$ and $q_b^*(x)$ are the auto and bus demand at location $x$ under the equilibrium conditions, respectively. $R^*$ represents the ratio of auto demand to total demand, and it is determined under equilibrium conditions, as outlined in Appendix A.

The travel demand of bus/auto at location $x$ is the cumulative demand density from the city boundary to location $x$, which is expressed as:

$$Q_j(x) = \int_x^A q_j(w)\, dw, \ \forall\, x \in (0, A]\,, \forall\, j = a, b, \tag{3}$$

where $j = b$ represents bus, and $j = a$ represents auto traffic.

Let $v(x)$ denote the total traffic volume including autos and buses at location $x \in (0, A]$:

$$v(x) = \frac{Q_a(x)}{O_a} + KF, \forall\, x \in (0, A], \tag{4}$$

where $K$ is the physical equivalent factor to convert the traffic effect of a bus into that of a few autos, and $F$ is the bus service frequency along the traffic corridor. $O_a$ is the average occupancy per auto, which can be obtained with Eq. (B6) in Appendix B.

### *2.2 System cost minimization model*

This section presents the objective function, which minimizes the total system cost by optimizing the bus frequency. The total system cost consists of both user costs and operator costs, as detailed in the following subsections.

#### *2.2.1 User cost*

User cost includes travel time, waiting time, passenger discomfort, and monetary expenses for both auto users and bus passengers, reflecting the total disutility experienced during trips along the corridor.

*When the traffic corridor is operating under the MTP*, the travel time per unit of distance at location $x$, denoted as $t_m^j(x)$, is obtained with the Bureau of Public Roads (BPR) function (Small, 1983). The function can be expressed as:

$$t_m^j(x) = t_0^j\left(1 + \alpha^j\left(\frac{v(x)}{n_{lane} * C}\right)^{\beta^j}\right), \ \forall\, x \in (0, A]\,, \forall\, j = a, b, \tag{5}$$

where $t_m^j(x)$ is the travel time per unit of distance for mode j under the MTP at location *x*, $t_0^j$ is the free flow auto/bus travel time per unit of distance, $C$ is the lane capacity, $n_{lane}$ is the number of lanes in the corridor, $\alpha^j$ and $\beta^j$ are the parameters under mixed traffic operation.

The travel time by bus/auto from location *x* to CBD $T_m(x)$ is:

$$T_m^j(x) = \int_0^x t_m^j(w)dw \ , \forall\, x \in (0, A], \forall\, j = a, b. \tag{6}$$

*When the corridor is operating under EBLP*, the auto and bus travel times per unit of distance are different, depending on the distance *x* to CBD. In contrast, the bus travel time per unit of distance is uniformly distributed along the corridor after the bus operating frequency *F* is determined. This can be expressed as:

$$t_e^a(x) = t_0^a\left(1 + \alpha^a\left(\frac{Q_a(x)}{O_a(n_{lane} - 1)C}\right)^{\beta^a}\right), \ \forall\, x \in (0, A]. \tag{7}$$

$$t_e^b(F) = t_0^b\left(1 + \alpha^b\left(\frac{KF}{C}\right)^{\beta^b}\right), \tag{8}$$

The travel time by bus from location *x* to CBD $T_e^b(x)$ under EBLP is:

$$T_e^b(x) = \int_0^x t_e^b(F)dw\,, \forall\, x \ \in (0, A], \tag{9}$$

where $t_e^b(F)$ is the bus travel time per unit of distance in the corridor, which is determined by the bus operating frequency $F$.

The travel time by auto from location *x* to CBD $T_e^a(x)$ under EBLP is:

$$T_e^a(x) = \int_0^x t_e^a(w)\, dw, \forall\, x\ \in (0, A], \tag{10}$$

*When the corridor operates under HOVLP*, both the high-occupancy vehicles and buses share the high-occupancy lane, resulting in the same travel time per distance unit. The average occupancy of low-occupancy autos is $O_{la}$ while the average occupancy of high occupancy autos is $O_{ha}$. The fractions of low-occupancy auto travelers $\varrho_l$ and high-occupancy auto travelers $\varrho_h$, relative to the total number of auto travelers can be determined based on the provided parameters, as explained in Appendix B. The travel time per distance unit for low-occupancy autos is expressed in Eq. (11). Buses and high-occupancy autos share the HOVL, resulting in identical travel times per unit distance, as specified in Eq. (12).

$$t_h^{la}(x) = t_0^a\left(1 + \alpha^a\left(\frac{\varrho_l Q_a(x)/O_{la}}{(n_{lane}-1)C}\right)^{\beta^a}\right),\ \forall\, x \in (0, A]. \tag{11}$$

$$t_h^{ha}(x) = t_h^b(x) = t_0^j\left(1 + \alpha^j\left(\frac{\varrho_h Q_a(x)/O_{ha}+KF}{C}\right)^{\beta^b}\right),\ \forall\, x \in (0, A], \forall\, j = a, b. \tag{12}$$

where $la$ and $ha$ in Eq. (11) and (12) represent low-occupancy auto and high-occupancy auto.

The travel time by bus and high-occupancy vehicle from location *x* to CBD $T_h^j(x)$ under HOVLP is:

$$T_h^b(x) = \int_0^x t_h^b(w) dw\ , \forall\, x \in (0, A], \tag{13}$$

The travel time by low-occupancy auto and high-occupancy auto from location *x* to CBD $T_h^j(x)$ under HOVLP can be represented by $T_h^{la}(x)$ and $T_h^{ha}(x)$, respectively:

$$T_h^\tau(x) = \int_0^x t_h^\tau(w)\, dw, \forall\, x\ \in (0, L], \forall\, \tau = la, ha. \tag{14}$$

The disutility of traveling by bus from location $x$ to CBD consists of the waiting time cost, the in-vehicle travel time cost, the in-vehicle passenger crowding discomfort cost, and the fare. It can be represented as:

$$U_p^b(x) = \varepsilon_w w_p^b(x) + \varepsilon_t^b T_p^b(x) + G_p^b(x) + f_b, \forall\, p = m, e, h, \forall\, x \in (0, A], \tag{15}$$

where $T_p^b(x)$ is the in-vehicle travel time from location $x$ to CBD, which can be determined with Eq. (6), (9), and (13). $\varepsilon_t^b$ is the value of in-vehicle time of bus travelers, $\varepsilon_w$ is the value of waiting time, and $f_b$ is the transit fare.

The average waiting time of passengers at location $x$ can be expressed by the volume-delay function (Li et al., 2009 and 2012):

$$w_p^b(x) = \frac{\gamma_1}{F_p} + \frac{\gamma_2}{F_p}\left(\frac{Q_b(x)}{O_b F_p}\right)^{\gamma_3}, \forall\, p = m, e, h,\ \forall\, x \in (0, A], \tag{16}$$

where $F_p$ is the bus frequency under policy $p$ (including MTP, EBLP, and HOVLP) in the corridor, and $\gamma_1$, $\gamma_2$, and $\gamma_3$ are positive calibrated parameters that depend on the distribution of bus headway and passenger arrival times.

Passenger discomfort $G_p^b(x)$ is generally affected by the degree of crowding at bus stops and onboard vehicles, which plays a significant role in mode choice because passengers may abandon the transit mode in favor of private auto travel when the discomfort becomes intolerable (Lam et al., 1999). The crowding discomfort cost per unit of in-vehicle travel time at location $x$ is expressed in monetary units and is assumed to increase with passenger demand at that location $x$, as shown by Eq. (17) (Huang et al., 2002):

$$g(x) = \iota_1\left(Q_b(x)\right)^2 + \iota_2\left(Q_b(x)\right), \forall\, x \in (0, A], \tag{17}$$

where $\iota_1$ and $\iota_2$ are positive parameters.

The total in-vehicle crowding discomfort cost $G_p^b(x)$ of a passenger traveling from location $x$ to CBD can be expressed as:

$$G_p^b(x) = \int_0^x g(w) t_p^b(w) dw \ , \forall\, p = m, e, h,\ \forall\, x \in (0, A], \tag{18}$$

where $t_p^b(w)$ is travel time per distance unit for buses under various policies, as defined by Eqs. (5), (8) and (12).

The total cost of the bus users $C_p^b$ under different policies can be obtained by accumulating the disutility from location $x$ to CBD for bus users:

$$C_p^b = \int_0^A U_p^b(x)\, q_b(x) dx, \tag{19}$$

The travel disutility of the auto mode is the sum of travel time cost and monetary cost (e.g., gas and parking cost). Therefore, under the MTP and EBLP, the travel disutility of the auto mode is formulated as Eq. (20):

$$U_p^a(x) = \varepsilon_t^a T_p^a(x) + \frac{\left(\delta_f^a + \delta_v^a x\right)}{O_a},\ \forall\ p = m, e\ , \forall\ x \in (0, A], \tag{20}$$

where $\delta_f^a$ and $\delta_v^a$ are the fixed cost (e.g., parking cost) and variable cost (e.g., gas cost per unit of distance), while $\varepsilon_t^a$ is the value of in-vehicle time of auto users.

Under the HOVLP policy, travel time and occupancy differ between the general lane and the HOV lane compared to the MTP. The travel disutility of the auto mode under HOVLP is expressed as:

$$U_h^a(x) = \varepsilon_t^a T_h^j(x) + \left(\delta_f^a + \delta_v^a x\right)/O_a^j, \forall\, j = la, ha, \forall\ x \in (0, A], \tag{21}$$

where $O_a^j$ is different auto occupancy in the general lane (low-occupancy vehicle lane) and HOV lane.

The total travel cost of auto users $C_p^a$ under different policies can be obtained by accumulating disutility from location $x$ to CBD for all auto users:

$$C_p^a = \begin{cases} \int_0^A U_p^a(x)\, q_a(x)dx, \ \forall\, p = m, e \\ \int_0^A U_h^{la}(x)\, \varepsilon_l q_a(x)dx + \int_0^A U_h^{ha}(x)\, \varepsilon_h q_a(x)dx, \ \forall p = h \end{cases} \tag{22}$$

*2.2.2 Operator cost*

The operator cost includes both the bus operating costs and the implementation costs of EBL/HOVL signals, which is the hourly cost per vehicle multiplied by fleet size (Sun et al., 2017; Li et al., 2011). The bus operating cost $C_p^O$ under different policies includes the fixed cost $C_O^f$ and the variable cost $C_O^v$ related to bus fleet size:

$$C_p^O = C_O^f + C_O^v\left(2T_p^b(A)F_p\right), \tag{23}$$

The fleet size required for operating the service is determined by $2T_p^b(A)F_p$, where $T_p^b(A)$ represents the one-way in-vehicle travel time along the corridor under policy $p$ and the length of corridor $A$, and $F_p$ denotes the service frequency. The round-trip time $2T_p^b(A)$ can be determined with Eq. (6). $C_O^f$ and $C_O^v$ are the fixed and variable cost for bus operation, in dollars per hour.

The EBL/HOVL signal implementation cost $C_p^s$ can be expressed as:

$$C_p^s = C_p^f + C_p^v A, \forall\, p = e, h, \tag{24}$$

where $C_p^f$ and $C_p^v$ are the fixed and variable parameters for EBL/HOVL signal implementation cost.

*2.2.3 Total system cost*

The total system cost of the corridor operating under different policies ($p = m, e, h$) includes bus user cost $C_p^b$, auto user cost $C_p^a$, bus operating cost $C_p^O$, and signal implementation cost of EBL and HOVL $C_p^s$. Therefore, the system cost minimization model can be formulated as:

$$min\ C_{p_{\{F_p\}}} = \begin{cases} C_p^b + C_p^o + C_p^a, \ \forall\ p = m, \\ C_p^b + C_p^o + C_p^a + C_p^s, \ \forall\ p = e, h. \end{cases} \quad (25)$$

$$s.t. F_p \geq \frac{\int_0^A q_b(x)dx}{o_b}. \quad (26)$$

The capacity constraint is defined in Eq. (26), ensuring that the transit service supply satisfies travel demand. The first-order optimality conditions for the cost minimization model under different policies are outlined in Appendix C, and a solution algorithm is applied to solve this system cost minimization model. Given the demand density, the mode choice between bus and auto is determined once the equilibrium condition is achieved, as described in Appendix A. In this algorithm, $R$ represents the ratio of auto demand to total demand in Eq. (2), serving as the mode choice at equilibrium. The solution algorithm for the cost minimization model is provided below.

**Algorithm**. Compute the system cost with corresponding optimal mode choice proportion and bus frequency.

**Require**: travel demand density $q$, given mode choice proportion range $\boldsymbol{R}$, given bus frequency range $\boldsymbol{F}$

**Ensure:** the algorithm returns system cost of policy with mode choice proportions $R^*$ and optimized bus frequency $F^*$

```
Define function C(q, R, F)
   Compute function R(q, R, F)
     Define function R(q, R, F)
       for R ∈ R do
          Compute function F(q, R, F)
            Define function F(q, R, F)
              for F ∈ F do
                   Compute function C(q, R, F)
                 Until  C(q, R, F) is minimal
              end for
            return  F* = argminC(q, R, F)
            Compute function C(q, R, F*)
         until C* = C(q, R, F*) is minimal
   end for
```

**return** $R^* = argminC(q, R, F^*)$

**return** $C^*$

## 3. Dynamic policy scheduling model

EBL and HOVL policies have been implemented in many cities worldwide. However, neither policy is suitable for continuous application throughout the day, as travel demand fluctuates significantly by time period. To address this issue, we propose a dynamic scheduling model that optimizes the deployment of different lane management policies, namely, the Mixed Traffic Policy (MTP), Exclusive Bus Lane Policy (EBLP), and High-Occupancy Vehicle Lane Policy (HOVLP). A practical example of this approach can be found on the I-110 Harbor Freeway in Los Angeles, California (Haddad, 2021). During weekday peak hours, the leftmost lane is reserved as an HOV lane (HOV3+) to accommodate carpools, vanpools, and registered high-occupancy vehicles, thereby maximizing person-throughput. At the same time, bus priority is maintained to ensure reliable transit service. Outside of peak hours, the lane is flexibly managed, often operating as an exclusive bus lane or as part of a dynamically priced high-occupancy toll (HOT) lane system, depending on demand conditions.

In the proposed model, the optimal bus frequencies for two different policies, $p1$ and $p2$, are denoted as $F_{p1}^*$ and $F_{p2}^*$, respectively. These are obtained through the solution procedure described earlier. The corresponding system costs are given by $C_{p1}\left(F_{p1}^*({q_0}^*)\right)$ and $C_{p2}\left(F_{p2}^*({q_0}^*)\right)$ that ${q_0}^*$represents the optimized demand-related design variables.

**Definition 1** (Switching threshold)

A threshold is determined by comparing each pair of policies; therefore, two thresholds exist with three policies. To simplify our equations and improve clarity, $p1$ and $p2$ are used to represent any two policies being compared. Here, $p1$ and $p2$ can be any two different policies from the set of three policies. Therefore, the system cost of one policy $p1$ dominates that of another policy $p2$ in the demand threshold if and only if the following conditions hold:

$$C_{p1}\left(F_{p1}^{*}({q_0}^{*})\right) \geq C_{p2}\left(F_{p2}^{*}({q_0}^{*})\right), \forall\, p = m, e, h, \tag{27}$$

where ${q_0}^{*}$ is the demand threshold determining the policy entry and exit (i.e., start and stop). The optimized bus frequency $F_p^{*}$ can be obtained from the first-order optimality conditions of the total cost function in Eq. (25). In this condition, the policy $p2$ is chosen because its total cost is lower; otherwise, the policy $p1$ is preferred. The first-order optimality conditions of the above cost function are given in Appendix C.

**Definition 2** (Switching timing for policy entry and exit)

Note that the costs of different policies are functions of the demand density. This naturally raises the question of whether there exists a demand density level at which different policies yield the same system cost. Since demand density fluctuates stochastically over time, it is modeled as a non-stationary stochastic process. To adapt to this dynamic nature, we adopt the Ornstein–Uhlenbeck (O-U) process, a continuous-time stochastic model widely used to describe mean-reverting dynamics (Guo et al., 2018). The O-U process, derived from Brownian motion, assumes that after a disturbance, the process gradually returns to its long-term mean (Ibe, 2013). It is characterized by two components: a drift term, which pulls the trajectory back toward the mean, and a shock term, which introduces randomness. The strength of mean reversion increases as the process deviates further from the mean, while the magnitude of shocks depends on the volatility of the historical data being modeled. Due to its ability to capture both randomness and mean-reverting tendencies, the O-U process is widely applied in financial forecasting, economics, and econometrics.

The stochastic differential equation (SDE) for the OU process:

$$dq(t) = \upsilon(\varpi - q(t))dt + \sigma dW_t, \tag{28}$$

where $\upsilon(\varpi - q(t))$ represents the drift term, where $\upsilon$ is the speed of mean reversion, describing how quickly deviations from the mean converge back. The parameter $\varpi$ is the long-term mean growth rate, interpreted as the average level of air travel demand growth approached in the long run. The stochastic component, $dW_t = W_{t+1} - W_t$, follows a standard normal distribution and acts as the shock term. (Geursen et al,. 2023).

The continuous-time solution of the O-U process can be expressed as (Sødal et al., 2008):

$$q_0(t) = q_0(0)exp(-vt) + \varpi\big(1 - exp(-vt)\big) + \sigma \int_0^t exp\big(-v(t-s)\big)\, dW_s. \quad (29)$$

This effectively showcases the dynamic and time-sensitive characteristics of the demand density variability. After establishing the demand density thresholds between any two policies as outlined in Definition 1, demand density is simulated by the O-U process, as described in Eq. (29). Given the dynamic nature of demand density, which fluctuates over time, there are two critical timings within each scheduling period for the different policies: one where the demand density exceeds the threshold and another where it subsequently falls below the threshold. Given the policies $p1$ and $p2$ and their associated parameters, demand density becomes a trigger factor for starting policy $p2$ from $p1$ at $t^1$ and exiting policy $p2$ to $p1$ at $t^0$, which is shown in Figure 3. The cumulative cost of different policies in each timing pair can be obtained as:

$$W_p(\rho(t)) = W_p\big(t^o, t^\xi\big) = \int_{t^o}^{t^\xi} C_{p_{\{F_p^*\}}}\big(q_0(t)\big)dt\,, \forall\, p = m, e, h, \quad (30)$$

where $W_p$ is the cumulative cost in the scheduled period for one policy, $C_p$ represents the system cost associated with policy $p$, including MTP, EBLP, and HOVLP, which can be computed using Eq. (25).

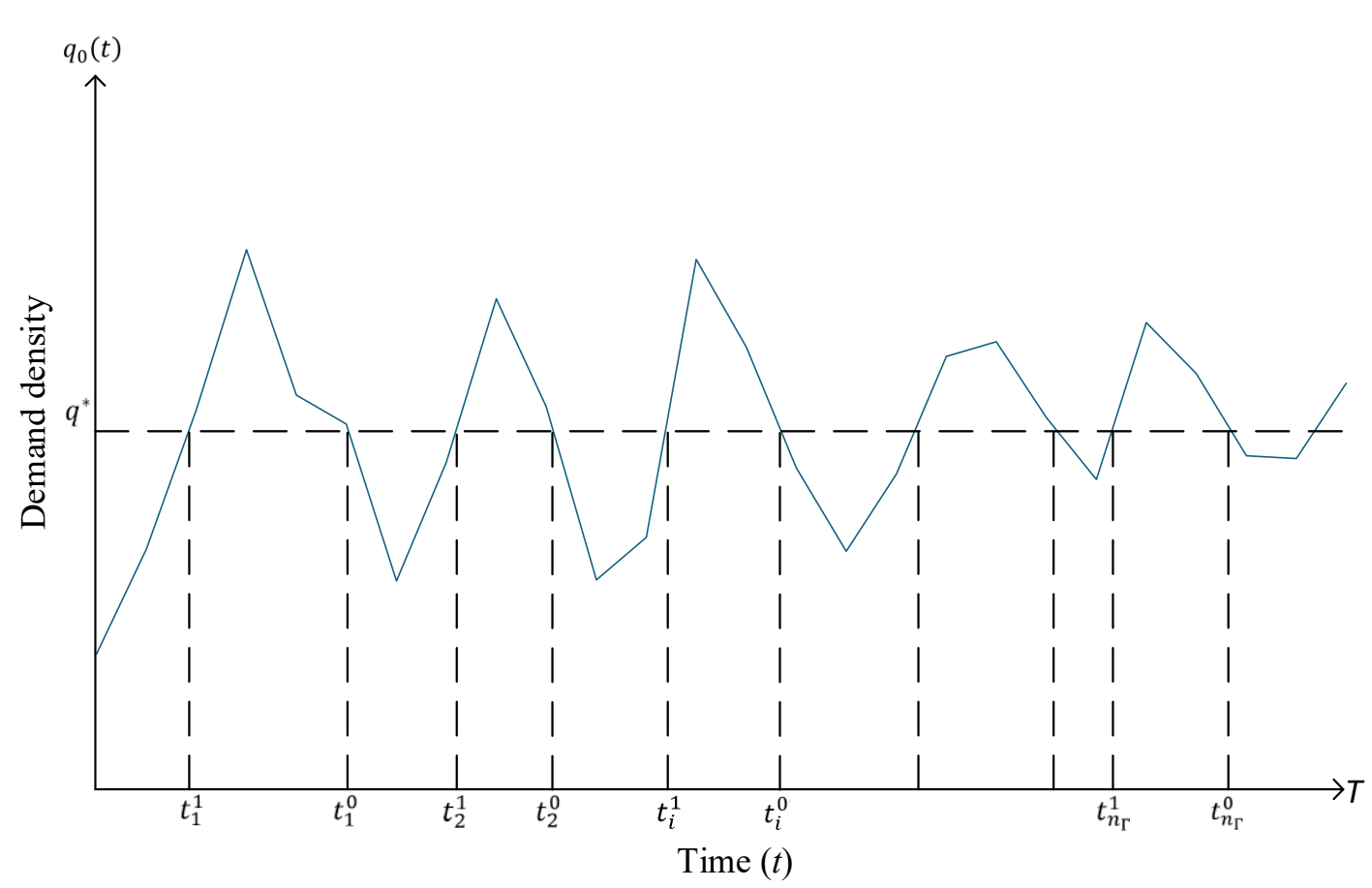

**Figure 3**. The timing of entry and exit.

The operating policy is defined by function $\rho(t) \in \{0,1\}$, where 0 represents an alternative policy with timing pair$(t^1, t^0)$ and 1 represents the current policy with timing pair $(t^0, t^1)$. At each decision point, the operator determines whether to switch from the previous policy state $\rho(t - dt) = 0$ (or 1) to the new state $\rho(t) = 1$ (or 0). The chosen policy directly affects the system cost $W_p(\rho(t))$ under the corresponding timing pair. An illustration of switching between policies is shown in Figure 3.

In the policy schedules problem, the objective is to identify feasible schedule periods that maximize total cost savings over the planning horizon $T$. Given two lane policies, $p1$ and $p2$, along with their associated parameters, the feasible lane policy schedule set is defined as $\Gamma = \{(t_1^1, t_1^0), (t_2^1, t_2^0), \cdots, (t_{n_\Gamma}^1, t_{n_\Gamma}^0)\}$ which must include all periods where the cost of one policy dominates the other, provided that Definition 1 is satisfied. Accordingly, the total cost savings can be formulated as:

$$Z = max_{\{\rho(t)\}} \sum_0^{n_\Gamma} \left(W_{p1}(\rho(t)) - W_{p2}(\rho(t))\right). \tag{31}$$

## 4. Numerical analysis

The data in Table 2 are used in the following numerical analyses to demonstrate the model.

**Table 2**. Data baseline

| Parameters | Definitions | Baseline values | Reference |
|---|---|---|---|
| $A$ | Distance between CBD and city boundary | 30 mi. | |
| $C$ | Lane capacity | 1500 veh/ hr | Basso et al., 2011 |
| $C_O^f, C_e^f, C_h^f$ | Fixed cost of bus/EBL/HOVL operation | \$300, 100, 500 | |
| $C_O^v, C_e^v, C_h^v$ | Variable cost of bus/EBL/HOVL operation | 20, 5, 10\$/hr | |
| $K$ | Factor converting a bus into an equivalent number of autos | 3 | Yao et al., 2018 |
| $f_b$ | Flat fare paid per bus trip | 1 \$/per trip | Sun et al., 2017 |
| $n_{lane}$ | The number of lanes in the traffic corridor | 3 lanes | Yao et al., 2018 |
| $O_{la}, O_{ha}$ | Occupancy rate of low/high occupancy auto | 1, 3 pax/veh | |
| $O_b$ | Bus capacity | 70 pax/veh | Guo et al., 2023 |
| $t_0^a, t_0^b$ | Free-flow auto/bus travel time per unit of distance | 0.05/ 0.025 hr./mi | HCM, 2010 |
| $\alpha^a, \alpha^b, \beta^a, \beta^b$ | Corrected parameters of the BPR function | 0.15/ 0.15/ 4/ 4 | HCM, 2010 |
| $\varepsilon_t^b, \varepsilon_t^a$ | Value of in-vehicle time of bus/auto travelers | \$15/ \$20 | |
| $\delta_f^a, \delta_v^a$ | Fixed/Variable cost of auto usage | \$2/ 0.3\$/mi | |
| $\iota_1, \iota_2$ | Positive parameters in discomfort function | 1.0e-6/ 5.0e-3 | Li et al., 2012 |
| $\mu$ | The fraction of low-occupancy autos over total autos | 0.6 | |
| $\gamma_1, \gamma_2, \gamma_3$ | Positive calibrated parameters in waiting time function | 0.5/ 0.05/ 2 | Li et al., 2009; |

### *4.1 Static analysis*

Figure 4 reflects the effect of demand density on total system costs of different policies, including MTP, EBLP, and HOVLP. To be more specific, Figures 4a and 4b display the trend of demand density in 200-1200 pax/hr/mi and 1200-2200 pax/hr/mi, respectively. Under the low-demand density conditions, congestion effects in general-purpose lanes are minimal, due to traffic volume remaining well below road capacity. As a result, both MTP and EBLP produce very similar system costs when the demand density is below 600 pax/hr/mi. Additionally, when congestion effects are minimal, the impact of bus and auto shares on total costs is also limited, as travelers are not incentivized to switch modes when neither option has a significant travel cost advantage over the other. However, as demand density increases, greater congestion occurs in both general-purpose and special lanes. When the cost of a particular mode decreases, it becomes more attractive, promoting a significant shift in travelers' mode choices and resulting in a notable impact on overall system costs. Consequently, the system costs of different policies change significantly when demand density exceeds 600 pax/hr/mi.

From Figure 4a, it is evident that HOVLP emerges as the more favorable policy at relatively low demand densities (below 1086 pax/hr/mi) since it has the lowest system cost. At low density, the HOVL remains largely uncongested, making it an efficient choice for low-demand periods where congestion effects are not a primary concern. As demand density increases (e.g., in Figure 4b), HOVLP becomes less favorable due to intensified traffic congestion, as evidenced by the sharp rise in auto travel cost shown in Figure 5b, which harms HOVL users. Since HOVL accommodates high-occupancy vehicles (HOVs) and buses, more passengers are affected compared with MTP and EBLP, thereby significantly increasing total system cost. Then, MTP becomes the more beneficial mode between 1086 pax/hr/mi and 2047 pax/hr/mi as HOVs enter the general-purpose lanes if the HOVL is more congested. Moreover, MTP proves beneficial by addressing the inefficiency of lane underutilization that can arise under dedicated-lane policies. Congestion is evenly distributed across all lanes under moderate demand, reducing the overall system cost by preventing bottlenecks specifically in HOVL or EBL. When demand density is high in Figure 4b, bus travelers tend to increase as potential travelers seek alternatives to avoid delays in general-purpose lanes. EBL enables buses to bypass congestion, reducing system costs by ensuring rapid public transit, which is essential in high-demand periods when road capacity is strained. Therefore, EBLP emerges as the most favorable policy in scenarios with extremely high demand density (above 2047 pax/hr/mi). Under MTP, where all types of vehicles share the lanes, congestion impacts all vehicles equally as demand grows. However, due to the higher occupancy rates of HOVs and buses, the system cost under HOVLP increases more rapidly than under MTP. As a result, the intersection point between HOVLP and EBLP occurs earlier than the intersection between MTP and EBLP. In summary, it is preferable to use HOVLP at low demand densities, switch to MTP as demand increases, and ultimately shift to EBLP during periods of very high demand.

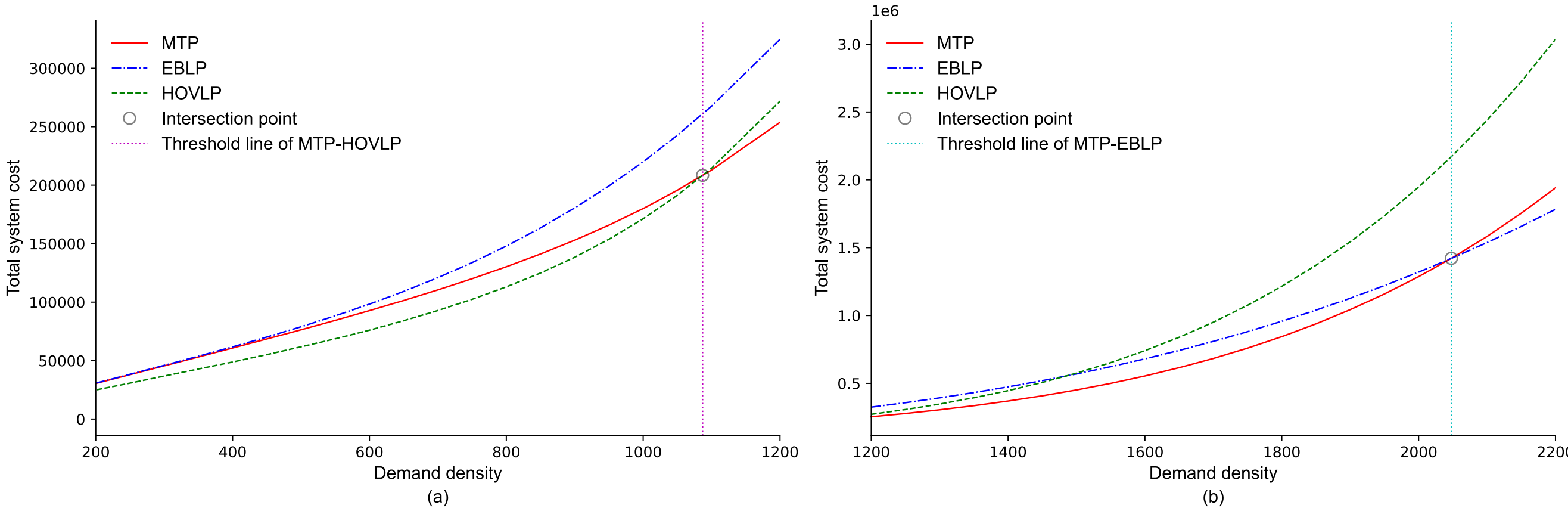

**Figure 4**. Total cost of different policies.

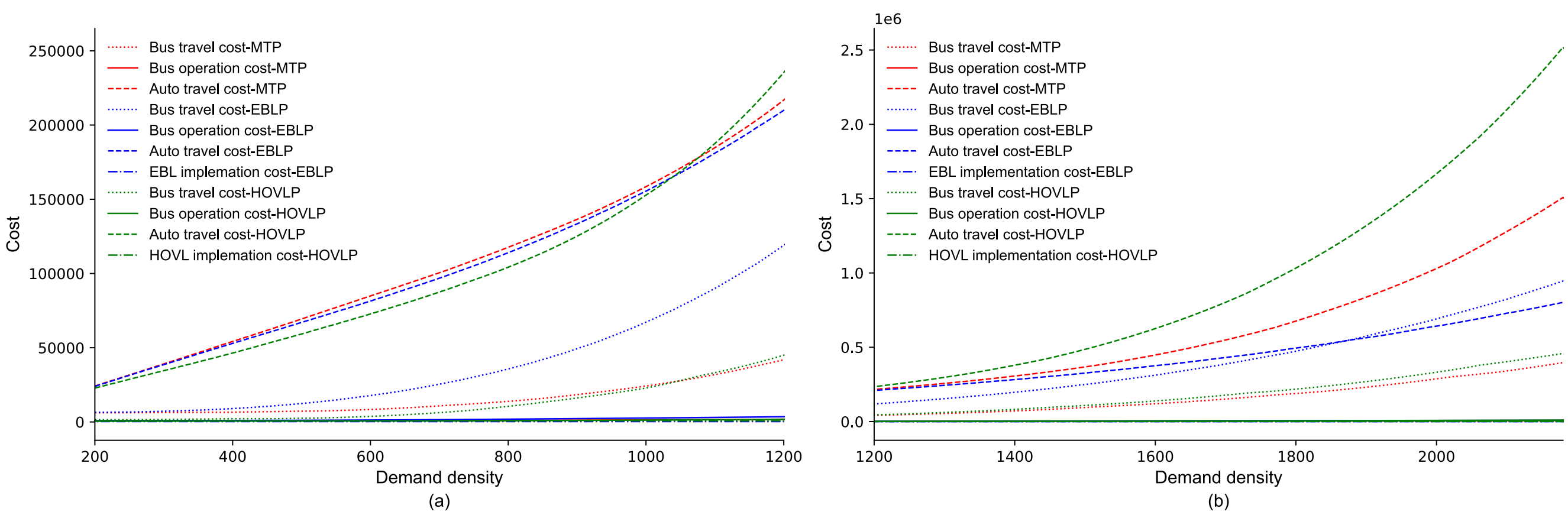

**Figure 5**. Breakdown of cost components.

Figure 5 shows the trend of different cost components as demand density rises. For a more detailed view, two graphs (Figure 5a and b) are presented from low to high demand density range. The mode split between bus and auto obtained through the user equilibrium mode choice model is shown in Figure 6. Comparing Figures 4 and 5 indicates that the total system cost is primarily affected by the auto travel cost (dashed line) and bus travel cost (dotted line). According to Figures 5a and b, bus operating cost (solid line) and lane implementation cost (dash dot line) have the least impact on the total cost compared to the travel cost of auto and bus, and it is difficult to observe their difference due to the magnitude difference from the other two costs. In Figure 5a, the auto travel cost dominates because most travelers opt for auto travel, as shown in auto share in Figure 6. As demand density increases from Figures 5a to 5b, a larger share of travelers shifts to bus travel, as indicated by the increased bus share under EBLP in Figure 6. This occurs because, as demand density increases, buses experience less congestion under EBLP, even as congestion intensifies in the general-purpose lanes. When demand density is high, bus travel cost rises due to more bus users. Meanwhile, auto travel cost under EBLP (blue dashed line) becomes lower than with other policies because the increase in bus demand decreases the total cost of auto travelers. In Figure 6, the auto share increases when the demand density is relatively low (below 400 pax/hr/mi) for all policies. Since the increased demand at low demand densities does not cause congestion in the corridor, the auto share continues to grow. The highest bus share of EBLP indicates that the implementation of EBL attracts more bus users. As the demand density increases, congestion occurs in the general-purpose lanes, making the advantage of EBLP more apparent and thereby increasing the bus share. The general trends of mode split are similar for MTP and HOVLP, indicating that the travel costs of autos and buses have the same trend as the total cost.

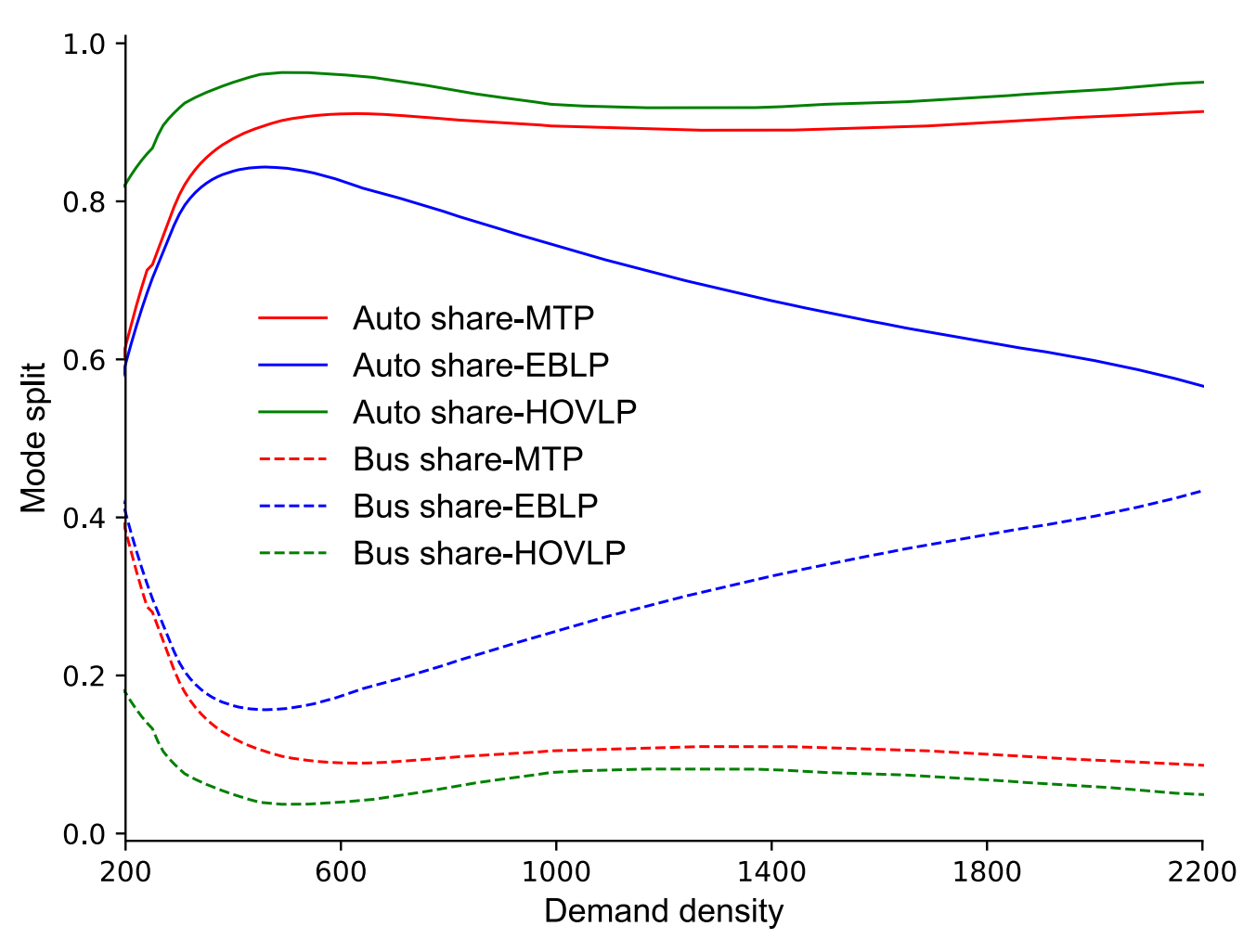

**Figure 6**. Mode split in the bi-modal traffic corridor.

### *4.2 Effect of the road capacity on total cost*

Figure 7 illustrates the variation in total system cost with demand density across different road capacities, ranging from 1200 to 1800 vehicles per hour. A larger $C$ indicates that more road resources can be utilized within the given time. As road capacity increases, the total system cost for all three policies decreases, as indicated by the progressively lower range of the Y-axis in Figures 7a1, 7b1, and 7c1. The cost difference between EBLP and MTP widens, and the intersection point changes from 1467 pax/hr./mi in Figure 7a2 to 2047 pax/hr./mi in Figure 7b2. This occurs because, in EBLP, an increase in road capacity delays the critical point at which the demand of auto travelers in general-purpose lanes exceeds the road capacity. As a result, EBLP becomes more advantageous, as the general-purpose lanes experience less congestion with increasing road capacity. For changes in the values of the Y-axis of the three lines, the decrease in total cost is most significant for EBLP as road capacity increases for all three policies. Additionally, Figures 7a1, 7b1, and 7c2 indicate that the cost difference between HOVLP and MTP widens before the intersection, and their intersection point changes from 864 pax/hr./mi. in Figure 7a1 to 1308 pax/hr./mi. in Figure 7c2. As road capacity increases, the point at which low-occupancy autos in general-purpose lanes under HOVLP reach the road capacity is delayed. Consequently, the benefits of HOVLP increase relative to MTP before their intersection point.

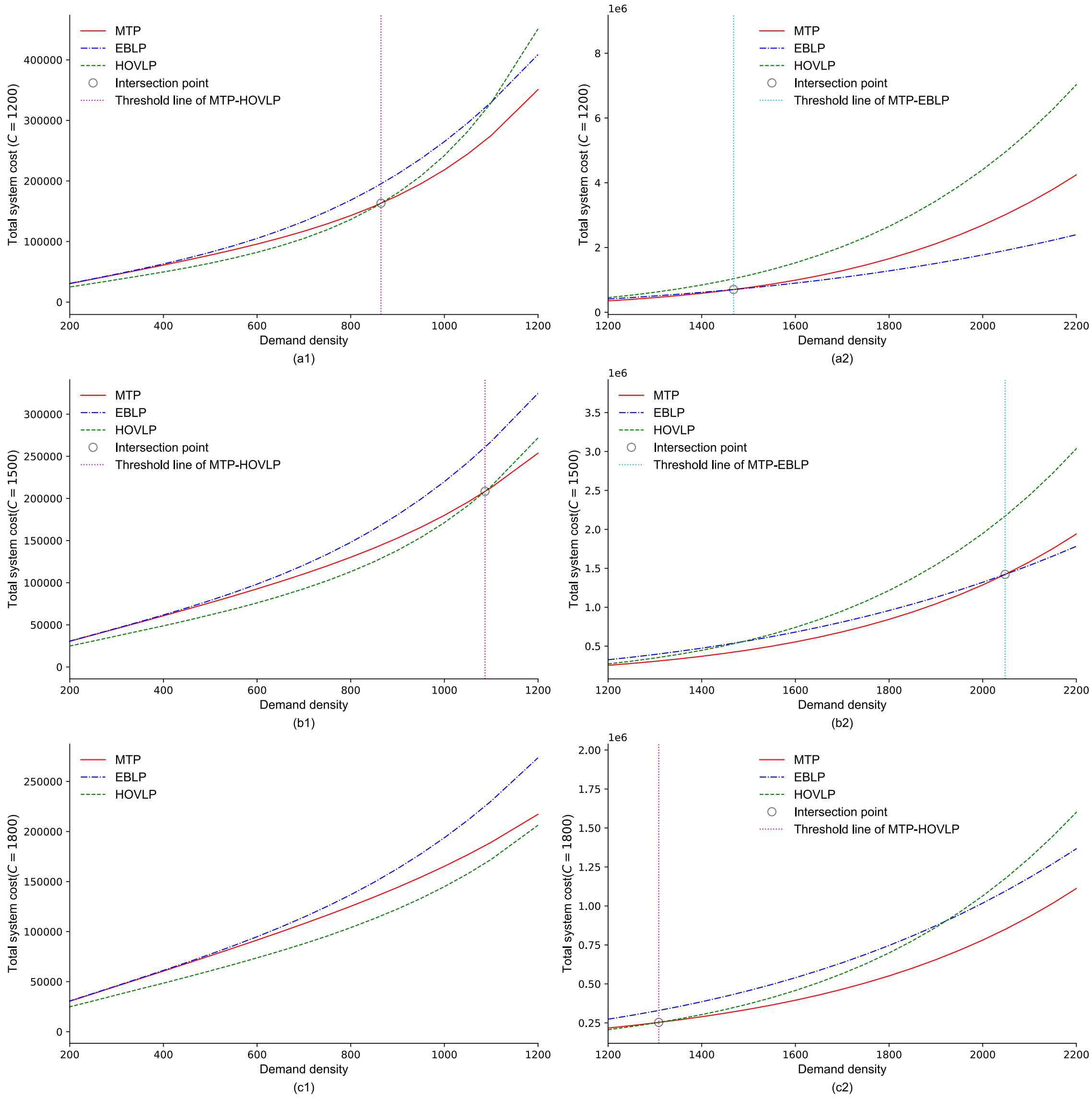

**Figure 7**. Impact of road capacity on total system cost under different policies.

### *4.3 Optimized frequency*

Figure 8 plots the optimized bus frequency for various policies. At first, the total cost initially decreases as bus frequency increases. This occurs because higher bus frequency reduces the waiting time for passengers, thereby lowering their travel cost and, consequently, the total cost. However, as bus frequency continues to increase, the operation costs of the buses also rise. Thus, at a later stage, the total cost begins to increase as bus frequency increases. From Figure 8, the optimized bus frequencies are 16 for MTP, 21 for EBLP and 11 for HOVLP, respectively. This difference arises because the extent to which bus is affected by auto varies across different policies. Under the HOVLP, lanes are shared between bus and HOVs, which might lead to more congestion in the HOVL than under the MTP because the travel cost is higher in Figure 9. This additional congestion could make it less efficient to operate buses at a higher frequency, as the benefits of increased frequency would be offset by delays caused by interactions with other vehicles in the lane. Therefore, the bus frequency is lower under the HOVLP than under the MTP. Under EBLP, the waiting time and travel time of bus passengers are not influenced by auto because buses operate separately in EBL. The optimized bus frequency is the highest compared with MTP and HOVLP. Therefore, bus operators need to offer higher frequencies when the reduction in travel cost exceeds the increase in operation cost. The predominant cost components that affect the optimized bus frequency are bus travel cost and bus operation cost, as shown in Figure 9, in which the left Y-axis represents the bus travel cost for different bus frequencies, while the right Y-axis indicates the bus operation cost.

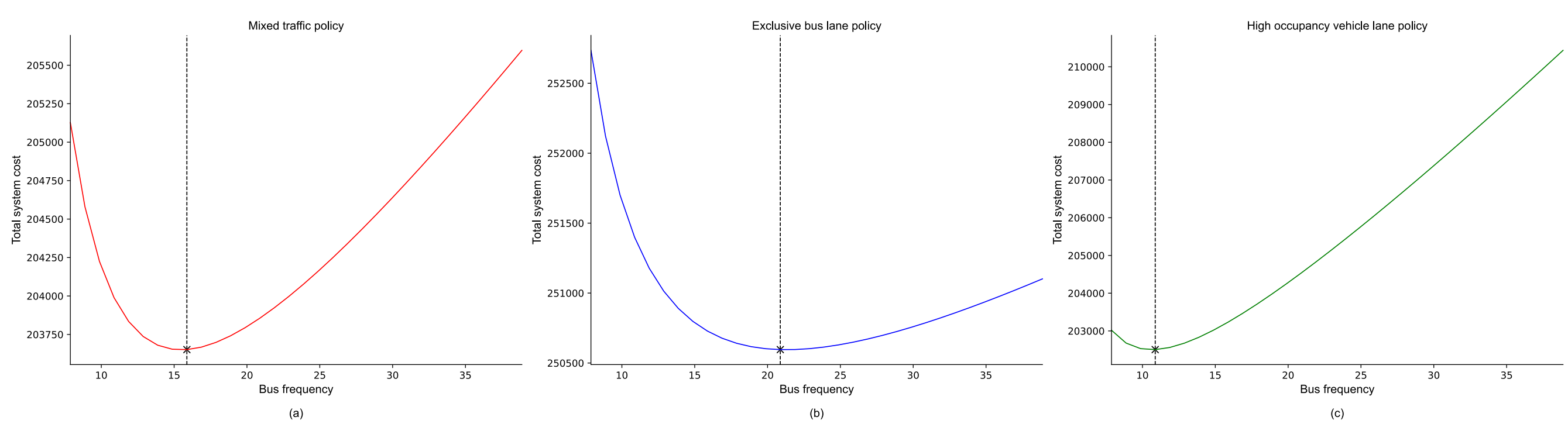

**Figure 8**. The optimized bus frequency under different policies.

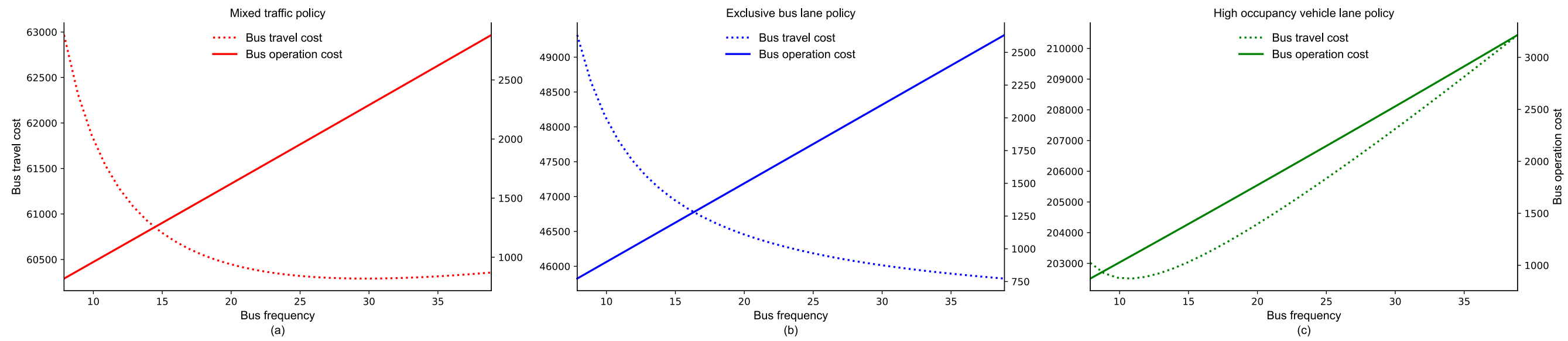

**Figure 9**. Optimized bus frequency, passenger travel cost, and operator cost under different policies.

### *4.4 Policy portfolio analysis*

In this section, we analyze three types of schedules using the proposed model, as summarized in Table 3. For Policy Schedules A and B, schedules are optimized under two policies (EBLP/HOVLP) with MTP. Two typical HOV and EBL corridors, I-5 and SR-99 in Seattle, are selected as case studies and simulated using the O-U process, with parameters calibrated from historical demand density data (Census Data 2020, Seattle). The schedules are evaluated by comparing policy performance and cost savings at different demand density thresholds. Furthermore, the combined policy schedule C is introduced, which incorporates transitions across the three policies. Model validation is conducted using simulated demand density, and Figure 10 provides a schematic illustration of the three schedules.

**Table 3**. Definition of three policy schedules.

| Proposed policy schedule | Features | Numerical experiment | Real case |
|---|---|---|---|
| Policy Schedule A | Switching between HOVLP and MTP | -- | Seattle I-5 corridor |
| Policy Schedule B | Switching between EBLP and MTP | -- | Seattle SR-99 corridor |
| Policy Schedule C | Switching among HOVLP, MTP, and EBLP | Simulated uncertainty demand density with O-U process | -- |

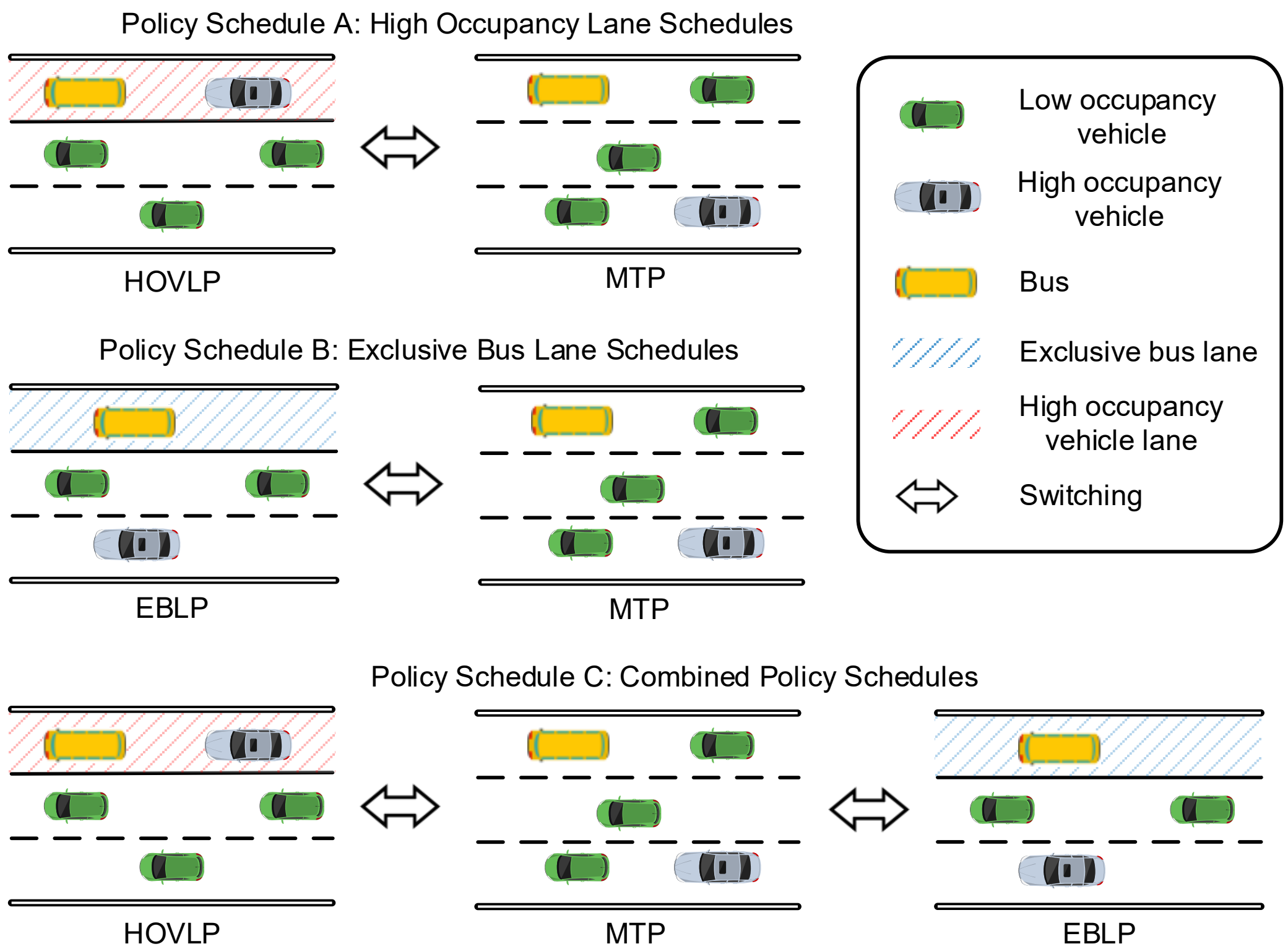

**Figure 10**. Schematic illustration of the three schedules.

### *4.4.1 Numerical experiment*

As noted previously, the demand density in the CBD over time is simulated using the O-U process (Sødal et al., 2008; Guo et al., 2017). Ten simulation trajectories CBD demand density, generated using Eq. (29), are presented in Figure 11. Based on these simulations, the total system costs of different policies can be evaluated. In the experiment, one trajectory from Figure 11 is selected, and the corresponding total costs of the three policies at different time points are illustrated in Figure 12. As demand density fluctuates, the relative performance of the three policies changes accordingly. Figure 12 highlights the most favorable policy over time and identifies potential switching points, represented by black circles, where the cost of one policy surpasses another. With real-time lane control, policies can be dynamically switched to minimize total cost. The vertical black lines in Figure 12 denote these switching times, derived under the infinite horizon setting using Eq. (31). For instance, the analysis shows that system costs vary with demand density, and the optimal policy differs across periods. Specifically, the HOVLP policy performs best under low-demand conditions, which occur only during two periods (7:33–8:02 and 14:31–15:45). Accordingly, HOVLP is scheduled for these periods in the timetable.

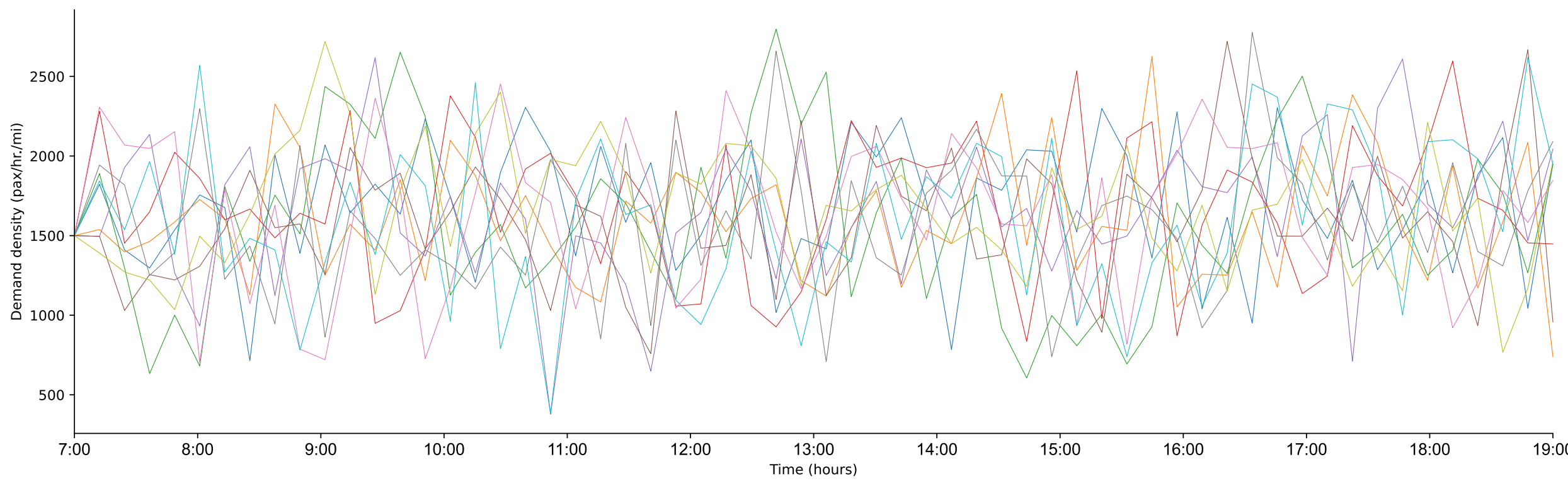

**Figure 11**. Ten demand density trajectories simulated by the O-U process.

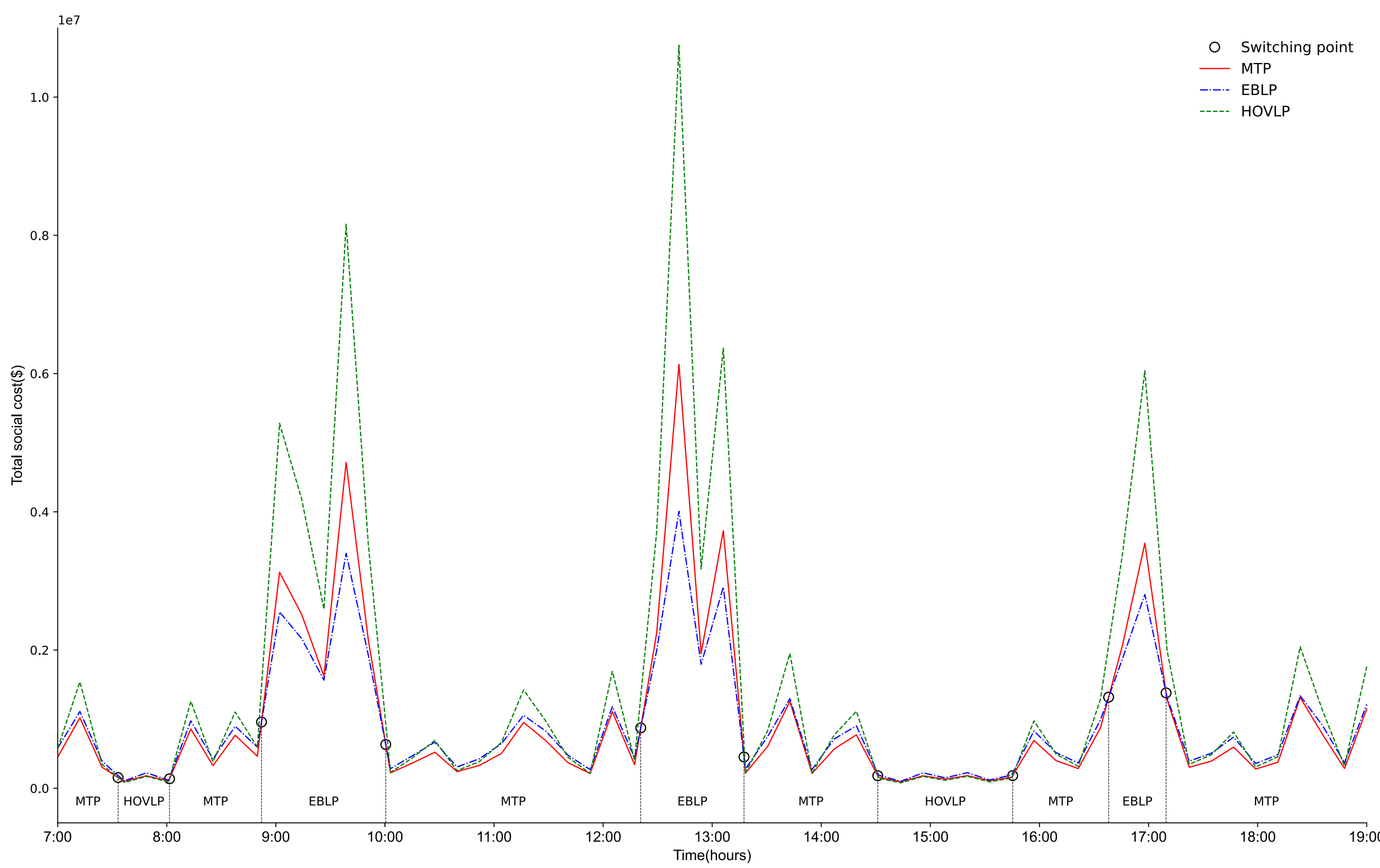

**Figure 12**. Optimal switching results among three policies.

The cumulative costs of different policies over the 12-hours from 7 am to 7 pm on workdays are summarized in Table 4, where we find that the cumulative cost of policy schedule C (combined policy schedules) decreases compared to the single policy (MTP, EBLP, and HOVLP only). Compared to single static policies (MTP, EBLP, and HOVLP only), the proposed combined policy reduces the cost by 11.9%, 6.7%, and 43.6%, respectively.

**Table 4.** Cumulative system cost of different schedules in the simulated policies.

| **Policy Schedule C** | **System cost under a single policy** (cost savings %) | | |
|---|---|---|---|
| Dynamic switching among MTP, EBLP, and HOVLP | MTP-only | EBLP-only | HOVLP-only |
| 9873593 | 11212514 (11.9%) | 10583157 (6.7%) | 17517774 (43.6%) |

### *4.4.2 Real-world case study*

To implement the proposed policy schedules model for Schedule A and Schedule B, Interstate 5 (I-5) and State Route 99 (SR-99) in Seattle are used as two real-world case studies. As shown in Figure 13, I-5 in the southern region of Seattle functions as a high-occupancy vehicle (HOV) lane corridor connecting suburban areas (Fife Heights) to the CBD, with a total corridor length of 27.7 miles. In contrast, SR-99 extends from northern Seattle to the CBD and has been designated as an exclusive bus lane (EBL) corridor, with a total length of 26.9 miles.

To apply our dynamic policy scheduling framework, we first estimate the parameters of the O-U process using observed demand data. Based on traffic volume data from the UW StarLab (DRIVE Net – UW, 2025), the estimated O-U parameters $(\upsilon, \varpi, \sigma)$ for I-5 and SR-99 are (4078, 1621, 20798) and (3756, 1028, 15317). The estimation procedure for I-5 is provided in Appendix F as an illustrative example. According to census data (Census Data 2020, Seattle), the study area population is approximately 54,280, with an average household size of four people, yielding about 13,570 households along the I-5 corridor.

Assuming 365 annual CBD trips per household and applying a peak-hour factor of 0.1, the resulting CBD demand density from I-5 is estimated at roughly 1,476 passengers per hour. State Route 99 (SR-99), extending from northern Seattle to the CBD, serves as an exclusive bus lane corridor with a total length of 26.9 miles. Using the same assumptions and census parameters, the CBD demand density from SR-99 is estimated at approximately 1,245 passengers per hour. According to WSDOT, the posted speed limits for I-5 and SR-99 are 60 mph and 35 mph, respectively; these values are applied as average speeds, while all other parameters are drawn from Table 2.

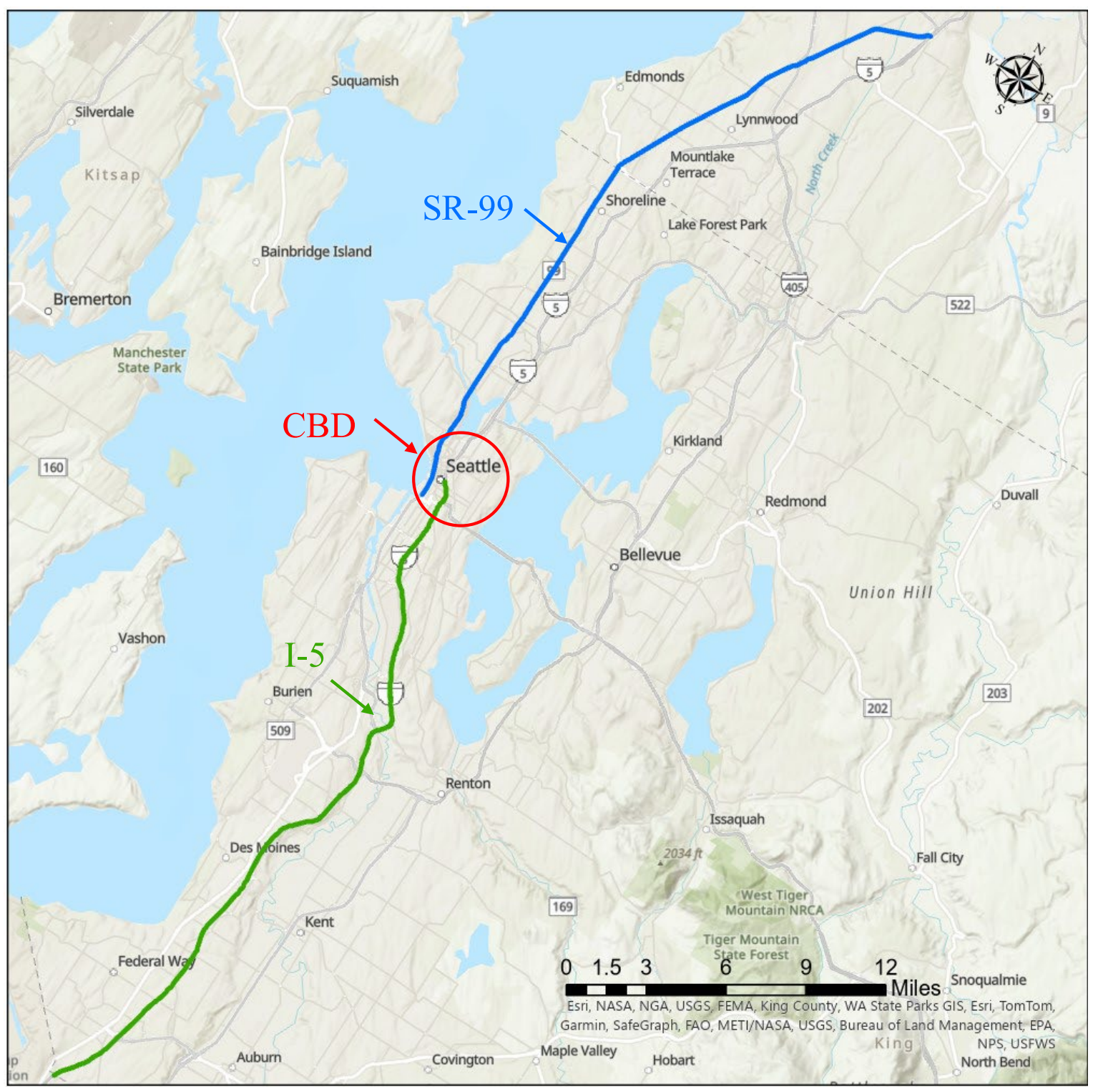

**Figure 13**. Two typical corridors in Seattle (USA Freeway System, 2025).

Using the same computation process as in the experiments, costs associated with different policies under Schedule A and Schedule B are computed to determine the threshold for switching. The fluctuations in demand densities for these policies are simulated using the O-U process based on the actual demand density of two corridors. Utilizing our proposed policy schedules model, schedules for these two policies are derived and presented in Table 5. The entry timing and exit timing are still defined for the switching timing between two different policies and the more beneficial policy option for each timing pair is detailed in Table 5. In real-world corridors, the conditions for switching among the three policies in the experiment are infeasible. Unlike in our experiments, our real-world case study utilizes existing HOVL and EBL corridors and applies our dynamic schedules according to the actual variation of demand. This demonstrates that our dynamic scheduling model can be effectively applied to existing EBL and HOVL corridors. Applying Schedule A and B in the existing HOVL and EBL, the cost savings from implementing the proposed policy schedules are computed and shown in Table 6. It is evident that the proposed schedules yield the lowest costs across both cases. Schedule A reduces costs by 13.2% and 32.5% compared to the MTP and HOVLP, respectively. Similarly, Schedule B reduces costs by 8.8% and 28.6% compared to the MTP and EBLP, respectively.

**Table 5**. Policy schedules of real-world case studies.

| **Schedule A applied in I-5** | | |
|---|---|---|
| Entry timing ($t^{entry}$) | Exit timing ($t^{exit}$) | Policy option |
| 7:00 | 8:11 | MTP |
| 8:11 | 9:23 | HOVLP |
| 9:23 | 16:49 | MTP |
| 16:49 | 18:12 | HOVLP |
| 18:12 | 19:00 | MTP |
| **Schedule B applied in SR-99** | | |
| Entry timing ($t^{entry}$) | Exit timing ($t^{exit}$) | Policy option |
| 7:00 | 8:21 | MTP |
| 8:21 | 12:15 | EBLP |
| 12:15 | 13:09 | MTP |
| 13:09 | 19:00 | EBLP |

**Table 6**. Cumulative cost savings of different policies in two case studies.

| | **Traditional schedules** | **Proposed** | **The cost savings compared** |
|---|---|---|---|

| | | | schedule | with traditional schedules | |
|---|---|---|---|---|---|
| | **MTP-only** | **HOVLP-only** | **Schedule A** | **MTP-only** | **HOVLP-only** |
| Cumulative cost ($/day) | 10511404 | 13517774 | 9127872 | 13.2% | 32.5% |
| | **MTP-only** | **EBLP-only** | **Schedule B** | **MTP-only** | **EBLP-only** |
| | 9160733 | 12703948 | 8354728 | 8.8% | 28.6% |

## 5. Conclusion and future studies

This study investigates the optimal dynamic configurations of road lanes in a commute corridor. A series of cost models, for the MTP, EBLP, and HOVLP, is formulated to compare their capabilities in enhancing traffic efficiency in a static demand setting. Such a static model is then extended to allow for dynamic switching among the three policies, depending on the varying real-time demand. Extensive numerical experiments and real-case studies are conducted to demonstrate the superiority of the dynamic switching scheme.

The results yield several key insights into optimizing lane management policies. A primary finding is that an optimal bus frequency can be effectively determined for each distinct policy (MTP, EBLP, HOVLP) within a traffic corridor, with the optimization primarily balancing the trade-off between bus operational costs and bus passengers' travel costs, including waiting and discomfort. The research also demonstrated that total system costs are dynamically influenced by the demand split between auto and bus usage. Specifically, HOVLP emerged as the most cost-efficient policy at demand densities below 1086 passengers per hour per mile (pax/hr/mi). MTP became the preferred option for demand levels ranging from 1086 to 2047 pax/hr/mi, while EBLP was found to be most advantageous when demand density exceeded 2047 pax/hr/mi. Furthermore, the study indicated that as road capacity increases, the resulting decrease in auto travel costs under an EBLP leads to a more substantial reduction in total system cost compared to MTP and HOVLP. This effect also highlighted how the cost difference between HOVLP and MTP can widen due to increases in auto travel costs under certain capacity conditions. Critically, the proposed combined policy schedules significantly reduced cumulative costs compared to the continuous implementation of individual static policies, achieving reductions of 11.9% against MTP, 6.7% against EBLP, and a notable 43.6% against HOVLP. Real-world case applications of two such dynamic schedules, Policy Schedule A and Policy Schedule B, further validated these benefits, demonstrating cumulative cost decreases of approximately 32.5% and 28.6%, respectively, when compared to cumulative costs under existing HOVL and EBL policies.

While this study offers valuable contributions listed in Section 1, it is subject to several limitations that also pave the way for future investigations. Firstly, the demand density was modeled using the O-U process. Future work could explore and compare the outcomes using alternative stochastic models, such as Geometric Brownian Motion, to assess model sensitivity and robustness. Secondly, the omission of intersection timing and delay, stemming from an assumption of continuous space for service, represents another area for refinement, which we extend in Appendix D. Thirdly, the current model does not account for uncertainties in road capacity, which can be affected by disruptions such as flooding or adverse weather. Integrating such supply-side uncertainties could be a significant extension. Finally, the analysis employed a many-to-one travel demand pattern, which, while useful, does not fully capture the complexities of more realistic many-to-many demand patterns characterized by multiple origins and destinations along the corridor. As a future direction, the model is extended to a many-to-many travel pattern framework in Appendix E. Future research endeavors will aim to incorporate a more comprehensive model to better represent these variations, thereby further improving the accuracy and applicability of policy evaluations.

**Acknowledgments:**

The authors Jiaqing Lu and Qianwen Guo were partially supported by the National Science Foundation (Grant No. 2400153).

**Declaration of interest statement:**

The authors declared no potential conflicts of interest with respect to the research, authorship, and/or publication of this article.

## Appendix A: Deterministic continuum equilibrium mode choice

To allocate demand effectively in a bimodal traffic corridor, it is crucial to reach a state where each traveler opts for the mode with the lower cost compared to the alternative mode. Once no traveler can find a more cost-effective travel mode, the congested corridor will achieve user equilibrium. According to Liu et al. (2009) and Li & Wang (2018), the user equilibrium conditions can be expressed as follows:

$$q_a(x) > 0 \Rightarrow U^a(x) \leq U^b(x), \tag{A1}$$

$$q_b(x) > 0 \Rightarrow U^b(x) \leq U^a(x), \tag{A2}$$

The above generalized travel disutility for the bus mode and auto mode, i.e., $U^a(x)$, and $U^b(x)$, are given by Eqs. (8) and (14), respectively. This definition states that at equilibrium, the individual generalized travel disutility by a mode at any location should be lower or equal to another mode if the mode is used by commuters at that location.

These are proven since, given the demand density $q(x)$, $x \in [0, A]$, the user equilibrium mode choice $\left(R^* = \frac{q_a{}^*(x)}{q(x)}, 1 - R^* = \frac{q_b{}^*(x)}{q(x)}\right)$ at location $x$ is an optimal solution of the following minimization problem (Liu et al., 2009):

$$min \int_0^A U^{a*}(x) q_a(x) dx + \int_0^A U^{b^*}(x) q_b(x) dx, \tag{A3}$$

$$subject\ to\ \ q_a(x) + q_b(x) = q(x), \forall x \in [0, A] \tag{A4}$$

$$q_a(x) \geq 0, \forall x \in [0, A] \tag{A5}$$

$$q_b(x) \geq 0, \forall x \in [0, A] \tag{A6}$$

where $U^{a^*}(x)$ and $U^{b^*}(x)$ are optimal user disabilities in the user equilibrium conditions.

Using the functional analysis method for minimization, we obtain the optimality conditions of the problem (A3) as follows:

$$q_a^*(x)(U^{a*}(x) - \lambda(x)) = 0, \forall x \in [0, A] \tag{A7}$$

$$U^{a*}(x) - \lambda(x) \geq 0, \forall x \in [0, A] \tag{A8}$$

$$q_a^*(x) \geq 0, \ , \forall x \in [0, A] \tag{A9}$$

$$q_b^*(x)\left(U^{b^*}(x) - \lambda(x)\right) = 0, \forall x \in [0, A] \tag{A10}$$

$$U^{b^*}(x) - \lambda(x) \geq 0, \forall x \in [0, A] \tag{A11}$$

$$q_b^*(x) \geq 0, \forall x \in [0, A] \tag{A12}$$

The above conditions show that $\lambda(x)$ is the equilibrium generalized travel cost of commuters leaving home at location $x$. Thus, all commuters at location $x$ have the identical and minimal generalized travel cost at equilibrium regardless of the mode chosen there. Eqs. (A7)-(A9) and (A10)-(A12) state the choice of travelers of the auto mode and bus mode along the corridor, respectively. For a given bus frequency, the travel demand distribution along the corridor can be determined by the mode choice equilibrium formulation Eqs. (A7)-(A12).

The equilibrium conditions in Eqs. (A7)-(A12) are derived under the assumption of deterministic demand $q(x)$, where mode choice stabilizes once no traveler can reduce their generalized travel cost by unilaterally switching modes. However, when demand evolves stochastically over time, such as when represented by an O-U process $q(x,t)$, the equilibrium must be interpreted dynamically. In this context, $q(x,t)$ fluctuates around a long-term mean with random perturbations, and the user equilibrium cannot be treated as a static allocation. Instead, the equilibrium allocation of demand $\left(q_a^*(x,t), q_b^*(x,t)\right)$ must be recomputed at each time step (or at discrete decision epochs) by solving the minimization problem in Eq. (A3) with the current realization of $q(x,t)$. This approach captures the adaptive nature of travelers' mode choice responses to short-term demand variations while maintaining the long-run equilibrium properties implied by Eqs. (A7)-(A12). In operational terms, this means that the corridor continually rebalances as demand fluctuates, and $\lambda(x,t)$ reflects the instantaneous generalized travel cost perceived by travelers at each location and time.

Under stochastic and time-varying demand $q(x,t)$, the user equilibrium at each time $t$ can be formulated as the following dynamic optimization problem:

$$\left(q_a{}^*(x,t), q_b{}^*(x,t)\right) = argmin_{q_a(x,t)+q_b(x,t)=q(x,t)} \int_0^A U^a(x,t) q_a(x,t) dx + \int_0^A U^b(x,t) q_b(x,t) dx, \quad \text{(A13)}$$

This formulation implies that the equilibrium is recomputed at each time step based on the instantaneous realization of stochastic demand, while ensuring mode choice consistency with the generalized travel disutilities $U^a(x,t)$ and $U^b(x,t)$.

## Appendix B: Derivation of the average occupancy

The fractions of low-occupancy auto travelers and high-occupancy auto travelers, relative to the total number of auto travelers, as well as the average occupancy per auto, can be determined based on the provided parameters. We utilize two equations to compute the fraction of low-occupancy auto travelers over the total number of auto travelers, denoted as $\varrho_l$, and the fraction of high-occupancy auto travelers over the total number of auto travelers, denoted as $\varrho_h$,. The given parameters include the total number of auto travelers $Q_a$, the occupancy of low-occupancy autos $O_{la}$, the occupancy of high-occupancy autos $O_{ha}$, and the fraction of low-occupancy autos over total autos $\mu$.

The total number of auto travelers is the sum of low-occupancy auto travelers $Q_{la}$ and high-occupancy auto travelers $Q_{ha}$:

$$Q_{la} + Q_{ha} = Q_a. \tag{B1}$$

For both occupancy rates of high-occupancy auto and low-occupancy auto, dividing by their respective fractions of the total number of auto travelers yields the total number of auto travelers, which can be expressed as follows:

$$\frac{Q_{la}}{O_{la} * \mu} = \frac{Q_{ha}}{O_{ha}(1-\mu)}. \tag{B2}$$

where $\mu$ is the fraction of low-occupancy autos over total autos.

By combining the two equations, we obtain the number of low-occupancy auto travelers $Q_{la}$ and high-occupancy auto travelers $Q_{ha}$, as expressed below:

$$\begin{cases} Q_{la} + Q_{ha} = Q_a, \\ \frac{Q_{la}}{O_{la} * \mu} = \frac{Q_{ha}}{O_{ha}(1-\mu)}. \end{cases} \tag{B3}$$

Then, $Q_{la}$ and $Q_{ha}$ can be solved as follows:

$$\begin{cases} Q_{la} = \frac{\mu Q_a O_{la}}{\mu O_{la} + (1-\mu)\, O_{ha}}, \\ Q_{ha} = \frac{(1-\mu) Q_a\, O_{ha}}{\mu O_{la} + (1-\mu)\, O_{ha}}. \end{cases} \tag{B4}$$

The $\varrho_l$ and $\varrho_h$ are determined as follows:

$$\begin{cases} \varrho_l = \frac{Q_{la}}{Q_a} = \frac{\mu O_{la}}{\mu O_{la} + (1-\mu)\, O_{ha}}, \\ \varrho_h = \frac{Q_{ha}}{Q_a} = \frac{(1-\mu)\, O_{ha}}{\mu O_{la} + (1-\mu)\, O_{ha}}. \end{cases} \tag{B5}$$

After determining $\varrho_l$ and $\varrho_h$, the average occupancy of auto $O_a$ can be obtained as follows:

$$O_a = \frac{Q_a}{\frac{\varrho_l Q_a}{O_{la}} + \frac{\varrho_h Q_a}{O_{ha}}} = \frac{O_{la} O_{ha}}{O_{ha} \varrho_l + O_{la} \varrho_h}. \tag{B6}$$

## Appendix C Derivation of optimized frequency

To obtain the optimized bus frequency, a sensitivity analysis-based approach and numerical search is used based on the first-order conditions of social cost minimization model. The details of a sensitivity analysis-based approach and a numerical search approach can be found in Ying and Yang (2005), and in Li et al. (2012). Taking MTP as an example, the system cost minimization problem can be expressed as:

$$min\ C_M(F_m) = C_m^b + C_m^o + C_m^a = \int_0^A \left(\varepsilon_w w_m^b(x) + \varepsilon_t^b \left(\int_0^x t_m^b(w)dw\right) + \int_0^x g(w) t_m^b(w)dw + f_b\right) q_b(x)dx + C_O^f + C_O^v(2T_m^b(A)F) + \int_0^A \left(\varepsilon_t^a \left(\int_0^x t_m^a(w)dw\right) + \delta_f^a + \delta_v^a x\right) q_a(x)dx. \quad \text{(C1)}$$

The first-order condition of above equation concerning bus frequency $F_m$ can be obtained by setting $\frac{\partial C_M(F_m)}{\partial F_m} = 0$, as follows:

$$\frac{\partial C_M(F_m)}{\partial F_m} = \int_0^A \left(\varepsilon_w \frac{\partial w_m^b(x,F_m)}{\partial F_m} + \varepsilon_t^b \frac{\partial \int_0^x t_m^b(w,F_m)dw}{\partial F_m} + \int_0^x g(w) \frac{\partial t_m^b(w,F_m)}{\partial F_m} dw\right) q_b(x)dx + 2C_O^v \left(\int_0^A \frac{\partial t_m^b(x,F_m)}{\partial F_m} dx F_m + \int_0^A t_m^b(x)dx\right) + \int_0^A \varepsilon_t^a \frac{\partial \int_0^x t_m^a(w,F_m)dw}{\partial F_m} q_a(x)\, dx = 0. \quad \text{(C2)}$$

The first-order condition of total system cost under EBLP and HOVLP, concerning bus frequency $F_e$ and $F_h$ can be derived similarly:

$$\frac{\partial C_E(F_e)}{\partial F_e} = \int_0^A \left(\varepsilon_w \frac{\partial w_e^b(x,F_e)}{\partial F_e} + \varepsilon_t^b \frac{dt_e^b(F_e)x}{dF_e} + \int_0^x g(w) \frac{dt_e^b(F_e)w}{dF_e} dw\right) q_b(x)dx + 2C_O^v \left(\frac{dt_e^b(F_e)A}{dF_e} F_e + t_e^b(F_e)A\right) = 0. \quad \text{(C3)}$$

$$\frac{\partial C_H(F_h)}{\partial F_h} = \int_0^A \left(\varepsilon_w \frac{\partial w_h^b(x,F_h)}{\partial F_h} + \varepsilon_t^b \frac{\partial \int_0^x t_m^b(w,F_h)dw}{\partial F_h} + \int_0^x g(w) \frac{\partial t_h^b(w,F_h)}{\partial F_h} dw\right) q_b(x)dx + 2C_O^v \left(\int_0^A \frac{\partial t_h^b(x,F_h)}{\partial F_h} dx F_h + \int_0^A t_h^b(x)dx\right) + + \int_0^A \left(\varepsilon_t^a \frac{\partial \int_0^x t_h^{ha}(w,F_h)dw}{\partial F_h}\right) ((1-\mu)\, O_{ha}/(\mu O_{la} + (1-\mu)\, O_{ha}))q_a(x)dx = 0$$

## Appendix D Intersection delay

This appendix provides the detailed formulation of intersection delay modeling used in the analysis and outlines how future research can incorporate signalized intersection control into the proposed framework. The current model quantifies intersection delays along the corridor based on evenly spaced intersections and pre-timed signal settings, distinguishing between different traffic management policies. Building on this foundation, we propose future extensions that integrate adaptive signal control and transit signal priority to more accurately capture the interaction between intersection operations and corridor-level performance.

It is assumed that all intersections are evenly spaced along the corridor, as illustrated in Figure D1. Therefore, the location of intersection $i$ can be represented by the expression $l_i = \frac{A}{n_{inter}+1} * i$, where $A$ denotes the total length of the corridor and $n_{inter}$ represents the total number of intersections.

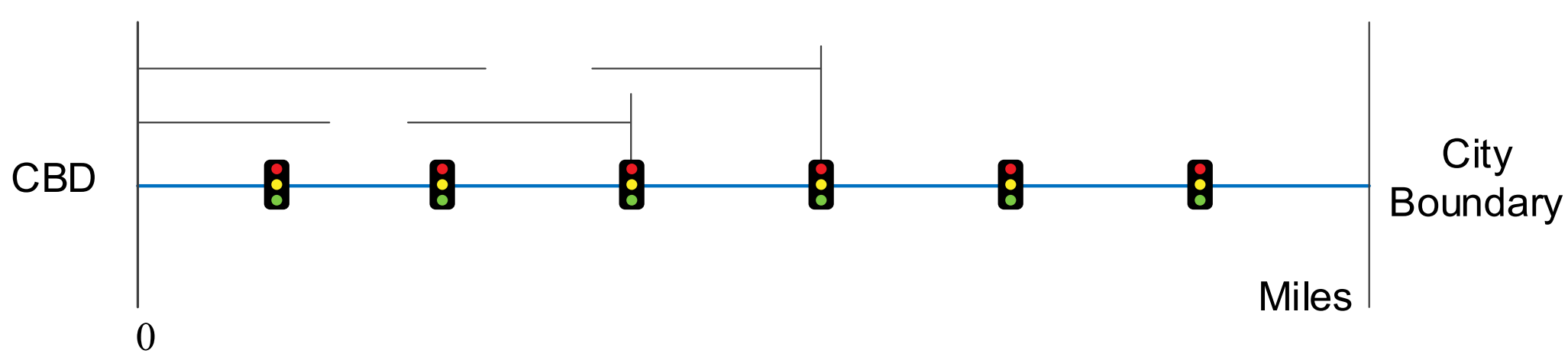

Figure D1. The distribution of intersections for the existing traffic corridor.

As stated in the delay equation (Dion et al., 2004; HCM 2010; Ge et al., 2013; Xie et al., 2023), the intersection delay per bus/auto at intersection $i$ under different policies is expressed as follows:

$$d_{pj}^{i} = \frac{C_{signal}(1-\lambda)^2}{2\left(1-min\left(1,\frac{\varpi_1}{\varpi_2}\right)\lambda\right)} + 900t_I\left(\left(\frac{\varpi_1}{\varpi_2}-1\right)+\sqrt{\left(\frac{\varpi_1}{\varpi_2}-1\right)^2+\left(8\kappa\varphi\frac{\varpi_1}{\varpi_2}\right)/(\varpi_2 t_I)}\right), \forall\, p = m, e, h, , \forall\, j = a, b,\ i \in [1, n_{inter}]. \quad \text{(D1)}$$

where $d_{pj}^{i}$ is the intersection delay per bus/auto at intersection $i$ under different policies, $C$ is the lane capacity, $C_{signal}$ is the cycle length, $\lambda$ is the green signal ratio of the lane, $t_I$ is the duration of analysis period, $\kappa$ refers to incremental delay factor accounting for pre-timed signal controller settings, and $\varphi$ is the adjustment factor for upstream filtering.

The intersection delay per bus and auto is identical at intersection $i$ due to the shared use of all lanes under MTP. The total volume at intersection $i$ is computed by summing the cumulative volume of autos and buses between intersection $i$ and $i+1$, which is expressed as $\int_{l_i}^{l_{i+1}} \frac{q_a(x)}{O_a} dx + KF/(n_{inter}+1)$. Under the EBLP and HOVLP, the delay functions differ due to variations in lane usage. In EBLP, buses and autos operate on separate lanes, whereas in HOVLP, high occupancy autos and buses use the dedicated lane, and low occupancy autos remain in the general-purpose lanes. Table D1 summarizes the intersection delay functions for both autos and buses under these distinct policies.

Table D1. Intersection delay functions under different policies.

| | Bus delay per vehicle ($d_{pb}^i$) | Auto delay per vehicle ($d_{pa}^i$) |
|---|---|---|
| MTP | $\varpi_1 = \int_{l_i}^{l_{i+1}} \frac{q_a(x)}{O_a} dx + KF/(n_{inter}+1)$; $\varpi_2 = n_{lane}C$ | |
| EBLP | $\varpi_1 = KF/(n_{inter}+1)$; $\varpi_2 = C$ | $\varpi_1 = \int_{l_i}^{l_{i+1}} \frac{q_a(x)}{O_a} dx$; $\varpi_2 = (n_{lane}-1)C$ |
| HOVLP | Bus/High occupancy auto | Low occupancy auto |
| | $\varpi_1 = (1-\mu)\int_{l_i}^{l_{i+1}} \frac{q_a(x)}{O_a} dx + KF/(n_{inter}+1)$;$\varpi_2 = C$ | $\varpi_1 = \mu\int_{l_i}^{l_{i+1}} \frac{q_a(x)}{O_a} dx$; $\varpi_2 = (n_{lane}-1)C$ |

The total intersection delay of bus travelers in the traffic corridor is determined by summing the delay times across all intersections. To compute the total delay time of bus travelers at intersection $i$, it is necessary to first determine the total number of bus travelers at that intersection. This is represented by the cumulative travel demand of bus travelers from the location of intersection $i$ to the city boundary. It can be expressed as follows:

$$D_p^b = \sum_{i=1}^{i=n_{inter}} \left( d_{pb}^i * \int_{l_i}^{A} Q_b(x) dx \right), \forall\, p = m, e, h, \quad \text{(D2)}$$

where $D_p^b$ denotes the intersection delay of buses under different policies.

Unlike Eq. (19), the total cost for bus users, $C_p^b$, under different policies must include the delay cost in the overall computation.

$$C_p^b = \int_0^A U_p^b(x)\, q_b(x) dx + \varepsilon_t^b D_p^b, \tag{D3}$$

## Appendix E Many-to-many demand framework

The study presented above simplifies passenger demand to a many-to-one framework, focusing on the travel demand from any location to and from the CBD at Node 0. As outlined in our future research, a key extension is the development of a many-to-many demand framework with modal competition between bus and auto, as illustrated in Figure E1.

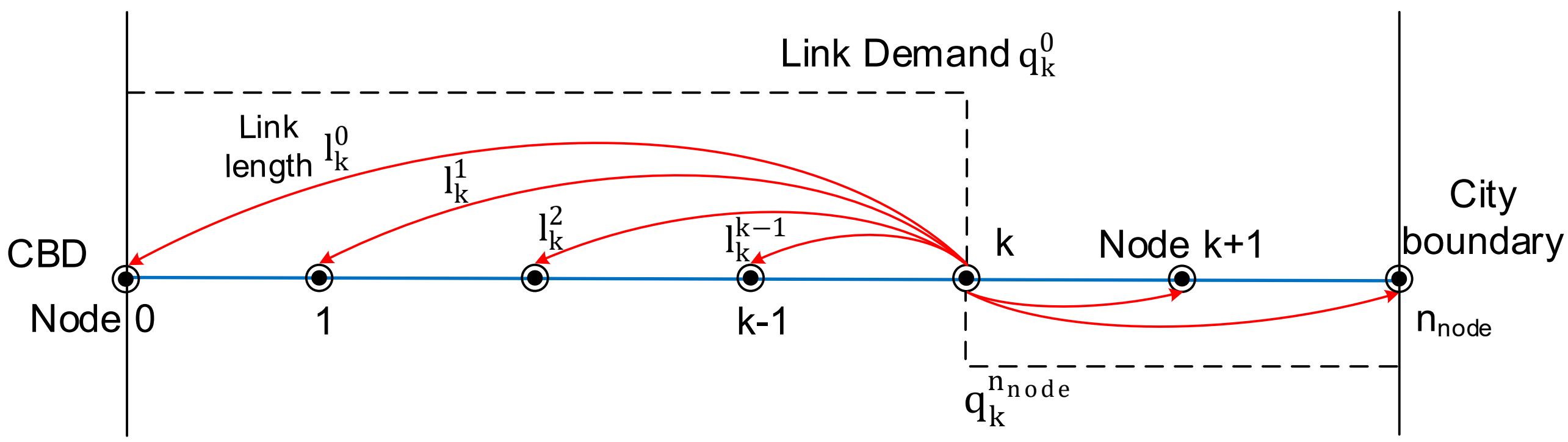

Figure E1. Illustration of many-to-many demand at a bimodal traffic corridor.

This extension would consider the total travel demand from node $i$ to node $j$, $q_i^j$, for all node pairs $(i, j)$ across the corridor, not just those demand from one location to the CBD. For each O-D pair, this total demand is split between bus users and auto users. The modal split for each O-D pair is determined by a user equilibrium (UE) condition for bus and auto, where no traveler can lower their travel cost by unilaterally switching modes. The UE solution can be obtained using the procedure described in Appendix A.

The total demand between any node $i$ and node $j$ can therefore be expressed as the sum of demand for each mode:

$$q_i^j = {q_a}_i^j + {q_b}_i^j. \quad \text{(E1)}$$

where ${q_a}_i^j / {q_b}_i^j$ is the portion of demand between nodes $i$ and $j$ traveling by auto/bus.

The model's objective function would still be to reflect the total passenger travel cost across both modes and all O-D pairs. The total travelers cost, $C_p^{a+b}$, is formulated as the sum of costs for all auto and bus users along the corridor:

$$C_p^{a+b} = \sum_{i=0}^{i=n_{node}} \sum_{j=0}^{j=n_{node}} U_p^a\left(l_i^j\right) {q_a}_i^j + U_p^b\left(l_i^j\right) {q_b}_i^j, \forall \ i, j \in [0, n_{node}]. \tag{E2}$$

where $n_{node}$ is the total number of nodes in the corridor, $l_i^j$ represents the generalized travel impedance (cost) between nodes $i$ and $j$, $U_p^a\left(l_i^j\right)$ denotes the travel cost for auto users, computed similarly to Eq. (20) based on the link length $l_i^j$ between nodes $i$ and $j$, and incorporating components such as travel time, fuel expenses, and parking fees, $U_p^b\left(l_i^j\right)$ denotes the travel cost for bus users, calculated analogously to Eq. (15) using the link length $l_i^j$, and including in-vehicle time, fares, waiting time, and discomfort costs associated with the location of discrete bus stop (node).

The system cost can then be integrated with other system costs (operation cost, policy implementation cost) which are determined by the formulas derived in the main body of this research. This future direction provides a more realistic and robust framework for traffic policy planning, directly addressing the need to evaluate service in a bimodal traffic corridor with complex, non-CBD-centric demand patterns.

**Appendix F Estimation of O-U process parameters**

The stochastic differential equation (SDE) for the O-U process is given by:

$$dq(t) = \upsilon(\varpi - q(t))dt + \sigma dW_t, \tag{F1}$$

To estimate the O-U process parameters from observed time series data, a least squares approach can be applied using the discrete-time solution of the SDE (Smith, 2010). At time $t + 1$, the process can be expressed as:

$$q_0(t+1) = exp(-\upsilon\Delta t)q_0(t) + \varpi\big(1 - exp(-\upsilon\Delta t)\big) - \sigma\sqrt{\frac{(1-exp(-2\upsilon\Delta t))}{2\upsilon}}\Delta W_t, \tag{F2}$$

where $\Delta t$ is sufficiently small, and $\Delta W_t$ represents independent and identically distributed Wiener process. This discrete-time formulation can be used to simulate long-term expected demand. $q_0(t = 0)$

Notably, Equation (F2) has the structure of a linear regression model, allowing straightforward parameter estimation via least squares (Smith, 2010):

$$q_0(t+1) = aq_0(t) + b + \epsilon_t. \tag{F3}$$

The parameters of the linear regression model are as follows:

$$\begin{cases} a = exp(-\upsilon\Delta t), \\ b = \varpi\big(1 - exp(-\upsilon\Delta t)\big), \\ sd(\epsilon_t) = \sigma\sqrt{\frac{(1-exp(-2\upsilon\Delta t))}{2\upsilon}}. \end{cases} \tag{F4}$$

where $sd(\epsilon_t)$ is the standard deviation of $\epsilon_t$.

By rearranging these equations, the parameters of the O-U process can be expressed as:

$$\begin{cases} \upsilon = \frac{\ln(a)}{\Delta t}, \\ \varpi = \frac{b}{1-a}, \\ \sigma = sd(\epsilon_t)\sqrt{\frac{2\upsilon}{(1-exp(-2\upsilon\Delta t))}}. \end{cases} \tag{F5}$$

We can fit this model to the historical traffic data, where the dependent variable is the traffic volume at the next time step $q_0(t+1)$ and the independent variable is the current traffic volume $q_0(t)$. The estimated regression coefficients, i.e., the intercept $a$ and slope $b$, along with the standard deviation of the $\epsilon_t$, are then used to compute the O-U parameters (Chaiyapo & Phewchean, 2017).

For the simulation, we extracted traffic volume data for all lanes of I-5 and SR-99 in the direction from the urban boundary toward the city center (the volume-increasing direction) from UW StarLab (DRIVE Net – UW, 2025) for the year 2016. To capture typical traffic demand dynamics, we focused on the 12-hour period between 7:00 and 19:00, which covers the morning and evening peak periods as well as the off-peak intervals. Figure F1 presents the results from 30 O-U process samples, estimated using the least squares approach, for the I-5 route in our case study. In Figure F1, the different gray lines represent the monthly average traffic volumes across the 12 months. The black line represents the annual average traffic volume, which serves as a smoothed baseline against which the monthly fluctuations can be compared. The initial values $q_0(0)$ in the O–U process are taken from the actual starting point of the historical data (red dashed line in Figure F1) to ensure that the simulations capture realistic demand variations. Therefore, simulated trajectories closely reproduce these observed variations, confirming that the O-U process provides a realistic stochastic representation of traffic demand dynamics on I-5.

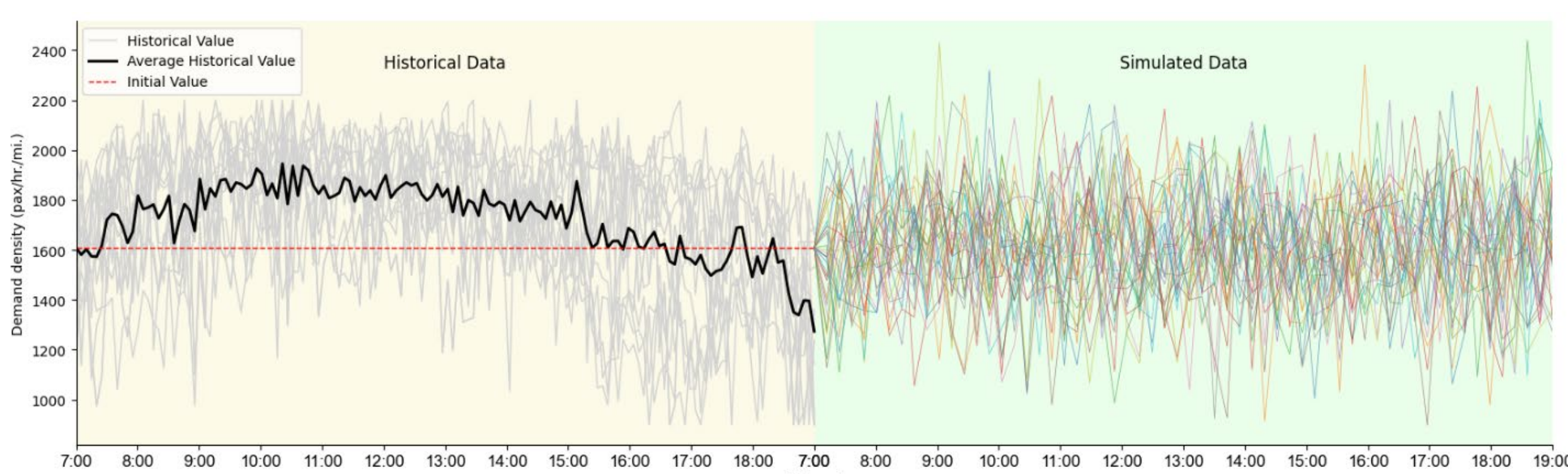

Figure F1.The case of I-5 demand simulation with O-U process.